\documentclass[11pt, a4paper, reqno]{amsart}

\usepackage[utf8]{inputenc}

\usepackage[letterpaper,hmargin=1.0in,vmargin=1.0in]{geometry}
\usepackage{amsmath,amsthm,amsfonts,amssymb,mathrsfs,stmaryrd,bigints}
\usepackage{dsfont}
\usepackage{mathtools}
\usepackage[usenames]{color}
\usepackage{MnSymbol}
\usepackage[colorlinks=true,linkcolor=blue]{hyperref}

\usepackage{tikz}
\usetikzlibrary{decorations.pathreplacing,calligraphy}

\usepackage{epstopdf}
\usepackage{amssymb}
\usepackage{setspace}
\usepackage{enumerate}
\usepackage{bigstrut}
\usepackage{multirow}
\usepackage{mathtools,xparse}
\usepackage{mathrsfs}
\usepackage{hyperref}
\usepackage{xcolor}
\usepackage [autostyle, english = american]{csquotes}
\newtheorem*{lemma*}{Lemma}
\newtheorem{theorem}{Theorem}[section]
\newtheorem{corollary}[theorem]{Corollary}
\newtheorem{lemma}[theorem]{Lemma}
\newtheorem{proposition}[theorem]{Proposition}
\theoremstyle{definition}

\newtheorem{remark}[theorem]{Remark}

\allowdisplaybreaks
\usepackage{chngcntr}

\newcommand{\E}{{\mathbb{E}}}

\newcommand{\ent}{h}
\newcommand{\Z}{\mathbb{Z}}

\newcommand{\R}{\mathbb{R}}

\newcommand{\e}{\varepsilon}
\newcommand{\de}\delta

\newcommand{\G}{\mathbb{G}}

\newcommand{\N}{\mathbb{N}}

	\renewcommand{\P}{\mathbb{P}}

\def\mc{\mathcal}
\newcommand{\cA}{\mc{A}}

\newcommand{\cE}{\mc E}

\newcommand{\co}{\scriptsize$\cO$\normalsize}

\renewcommand{\t}{\tilde}

\newcommand{\wt}{\widetilde}
\newcommand{\wh}{\widehat}

\newcommand{\la}{\leftarrow}

\DeclareMathOperator{\diam}{diam}

\DeclareMathOperator{\vol}{vol}

\newcommand{\f}{\frac}

\renewcommand{\setminus}{\backslash}
\newcommand{\sm}{\backslash}

\newtheorem{maintheorem}{Theorem}

\def\AA{\mc A}

\def\E{\mathbb E}

\def\Ld L
\def\L{\mathbb L}
\def\N{\mathbb N}
\def\PP{\mc P}

\def\P{\mathbb P}

\def\R{\mathbb R}

\def\G{\Gamma}

\def\Z{\mathbb Z}

\def\been{\begin{enumerate}}
\def\bee{\begin{example}}
\def\beit{\begin{itemize}}
\def\bel{\begin{lemma}}
\def\bepr{\begin{proposition}}
\def\bep{\begin{proof}}
\def\bet{\begin{theorem}}
\def\bec{\begin{corollary}}
\def\enc{\end{corollary}}

\def\bs{\boldsymbol}

\def\co{\colon}

\def\de\de 
\def\diam{\ms{diam}}

\def\dist{\ms{dist}}

\def\enen{\end{enumerate}}
\def\ene{\end{example}}
\def\enit{\end{itemize}}
\def\enl{\end{lemma}}
\def\enpr{\end{proposition}}
\def\enp{\end{proof}}
\def\ent{\end{theorem}}

\def\es{\varnothing}
\def\e{\varepsilon}
\def\ff{\infty}
\def\f{\frac}
\def\g{\gamma}

\def\la{\lambda}

\def\ms{\mathsf}
\def\Iso{\ms{Iso}}
\def\de{\delta}

\def\pa{\partial}

\def\sm{\setminus}
\def\su{\subseteq}
\def\s{\sigma}
\def\to \infty{\uparrow\infty}

\def\th{\theta}

\def\t{\tau}

\def\x{x}

\def\bef{\begin{figure}[!h]}
\def\enf{\end{figure}}

\def\ccf{\mathbf C^\ff}

\def\ccn{\mathbf C^n}
\def\cc{\mathbf C}

\def\vol{\ms{vol}}
\def\len{\ms{len}}
\def\poly{\ms{poly}}
\def\hull{\ms{hull}}

\def\hgp{\hat\g}

\def\intt{\ms{int}}
\def\wh{\widehat}
\def\wtpn{\Phi_n}

\def\aac{\AA_{n, \de}^{\ms{app}}}

\begin{document}

\title[Large deviations for  isoperimetry in 2D percolation]{Large deviations for the isoperimetric constant in 2D percolation}

\author{Christian Hirsch and Kyeongsik Nam}

\begin{abstract}
Isoperimetric profile describes the minimal boundary size of a set with a prescribed volume. 
Itai Benjamini conjectured that the isoperimetric profile of the giant component in supercritical percolation experiences an averaging effect and satisfies the law of large numbers. This conjecture was  settled by Biskup-Louidor-Procaccia-Rosenthal 
for 2D  percolation  \cite{blpr}, and later resolved  by Gold  for higher-dimensional lattices \cite{gold}.

However, more refined properties of the isoperimetric profile, such as  fluctuations and large deviations,  remain unknown. 
In this paper, we determine the large deviation probabilities of the isoperimetric constant in 2D supercritical percolation, answering the question in \cite[Open problem (6)]{blpr}. Interestingly, while the large deviation probability is of surface order in the entire upper tail regime, a phase transition occurs in  the lower tail regime, exhibiting both surface  and volume order large deviations. 
\end{abstract}

\address[Christian Hirsch]{Department of Mathematics, Aarhus University, Ny Munkegade, 119, 8000, Aarhus C, Denmark}
\address[Christian Hirsch]{DIGIT Center, Aarhus University, Finlandsgade 22, 8200 Aarhus N, Denmark}
\email{hirsch@math.au.dk}

\address[Kyeongsik Nam]{Department of Mathematical Sciences, KAIST, South Korea}

\email{ksnam@kaist.ac.kr}

 \subjclass[2010]{}


\maketitle

\tableofcontents

\section{Introduction}
\label{sec:intro}

The problem of isoperimetry concerns the smallest possible  perimeter of a set given its volume. While having its roots in metric geometry, isoperimetric inequalities have now a central place in various applications in diverse fields of mathematics. For instance, they appear in the spectral analysis of Laplacians, in the investigation of the mixing time of Markov chains and in heat kernel estimates \cite{cheeger,lk,var}.

While Markov chains on graphs  have been investigated for several decades, more recently there has been a vibrant study on  understanding stochastic processes on random networks. One prototypical example of such a random network is a (bond) percolation on graphs. The concept of isoperimetry has  played a crucial role in understanding the   phase separation phenomenon in percolation.  
Indeed, the macroscopic shape of  phase separation is described in terms of a solution to the particular isoperimetric
problem minimizing 
the surface free energy, called \emph{Wulff construction}. 
This was first justified in \cite{acc} for 2D percolation and then later established in higher-dimensional lattices \cite{bod,cerf1,cerf,cerf3,surface,surface1} (see \cite{dks,cerf0} for a survey on Wulff construction).

A major difficulty in understanding the isoperimetric profile on such random networks is that it fluctuates between different realizations in a complicated way. In \cite{bm}, the Cheeger constant (see \eqref{eq:blpr} for the definition) of a giant component in supercritical (bond) percolation on $B(n) = [-n,n]^d \cap \Z^d$ ($d \ge 2$) is shown to be of order $\Theta(\frac{1}{n})$ with high probability (see also \cite{remy,pete}). This was crucially used to deduce that the mixing time of the simple random walk on the giant component  is of order $\Theta(n^2),$ with the aid of   Lov\'asz-Kannan method \cite{lk}.

Itai Benjamini conjectured that the  (rescaled) isoperimetric profile of the giant component in supercritical percolation  satisfies the law of large numbers (LLN) in growing domains. The seminal result is obtained   by Biskup-Louidor-Procaccia-Rosenthal \cite{blpr} who resolved this  conjecture for the 2D lattice. Later, this  was fully settled by Gold  \cite{gold} for the higher-dimensional lattices $\Z^d$ ($d\ge 3$). We also refer to \cite{continuity}
 for a continuity property of  the isoperimetric profile in the percolation parameter.

While the  LLN for the isoperimetric profile in supercritical percolation is an essential milestone, it leaves open the question of more refined fluctuation and large deviation properties. In the context of the classical phase separation problem (Wulff construction) in the 2D Ising model, fluctuation properties of a droplet have been established  \cite{alexander,hammond3,hammond2,hammond1}. We also refer to \cite{basu} where the fluctuation  of a maximizing path is studied in the model of Poissonian last passage percolation which  reflects some features of isoperimetry.  However, fluctuation properties of  the isoperimetric profile in supercritical percolation are not well-understood  yet (see \cite[Section 1.4]{basu} for the conjectural fluctuation behavior).

Regarding large deviations, \cite[Open problem (6)]{blpr} proposes the question of finding a correct order of large deviation probabilities  for the rare event that the isoperimetric profile deviates from its typical value.
In particular, \cite[Open problem (6)]{blpr} tentatively hypothesize that similar to the situation in first-passage percolation (FPP), the upper and lower tail large deviation probabilities could exhibit different orders. Indeed, in the case of passage times in  classical 2D FPP problem (with bounded single-edge passage times), the upper and lower tail large deviation probabilities are of orders $e^{-\Theta(n^2)}$ and $e^{-\Theta(n)}$, respectively \cite{kestAsp}.

To the best of authors' knowledge, the only existing large deviation result regarding the  isoperimetric profile in percolation is the one obtained by Pete \cite[Theorem 1.2]{pete}: In supercritical percolation on $\Z^d$ with connectivity probability $p>p_c(\Z^d)$, denoting by $o$ and $\ccf$ the origin and (unique) infinite cluster respectively, there exist  constants $\alpha,c>0$ (depending on $p$ and $d$) such that 
\begin{align*}
    \P\Big( \inf \Big\{ \frac{|\partial_{\ccf} S|}{|S|^{1-1/d}} :  S \text{ connected},   o\in S \su \ccf,  n \le |S|<\infty \Big\} \le \alpha \Big) \le e^{-c n^{1-1/d}}
\end{align*}
(see \eqref{eq:blpr}  for the definition of edge boundary $\partial_{\ccf} S$).
In other words, it provides an upper bound on the ``much'' lower tail large deviation probability (i.e. the regime far from its typical value).
However, this result is not enough to entirely understand the  large deviation phenomenon: (a) The matching lower bound  is not provided.  (b) The lower tail large deviation is missing in the regime ``near'' the typical value.  (c) The upper tail large deviation is not examined.

In the present paper, we precisely examine the large deviation properties regarding the isoperimetric  constant of the giant component in 2D supercritical percolation. Indeed, we find the correct order of large deviation probabilities  in the entire rare-event regimes. 
The main contribution of our work is that in fact, large deviation probabilities are quite different from the expectation in \cite[Open problem (6)]{blpr}.  Precisely, we verify that the upper tail  large deviation probabilities are of surface order in the entire upper tail regime. Furthermore, in the lower tail regime, we explicitly identify the locations whose large deviation probabilities are of surface and volume orders respectively.  Although our arguments rely on the planarity of the  lattice, we expect that our surface  and volume order large deviation result extends to the higher-dimensional lattices as well (see Section \ref{high} for explanations).

  Our results build on surprising findings about the large deviations for the chemical distance in supercritical percolation \cite{gm}. 
   Chemical distance in  percolation  can be regarded as a FPP problem, but it violates the boundedness of the single-edge passage time. Indeed,   passage time of closed and open edges in percolation  can be considered as $\infty$ and 1, respectively. This is one of the key reasons why the dichotomy between upper and lower tails of  the passage time in FPP  \cite{kestAsp} does not extend to supercritical percolation. 
 {In addition, due to the intricate interplay of the volume and chemical distance in the isoperimetry problem, we observe a phase transition phenomenon in the lower tail for the isoperimetric profile, exhibiting both volume and surface order large deviations, unlike the FPP problem.}
 
    The rest of the manuscript is organized as follows. In Section \ref{sec:res}, we state Theorem \ref{thm:main} as our main result that concerns the large deviation properties of the isoperimetric constant in  2D supercritical percolation. We present some useful properties about percolation   and isoperimetric problems in Sections \ref{sec:pre} and \ref{sec:isop}. Then, in Sections \ref{sec:up000} and \ref{sec:low}, we prove the large deviations asymptotics in the upper and lower tails, respectively.  Finally, in Section \ref{sec:App}, we prove technical auxiliary results on the cluster density and chemical distances in supercritical percolation.
 
\section{Main results}
\label{sec:res}

Our work concerns supercritical Bernoulli bond percolation on $\Z^2$, where each of the edges in $\Z^2$ is retained independently with a probability $p \in ( 0,1)$. By the Harris-Kesten theorem \cite[Theorem 11.11]{grim}, for $p > p_c:= 1/2$ the resulting graph almost surely has a unique infinite connected component consisting of open edges, which is called $\ccf$. To lighten notation, we assume that $p \in  (1/2,1)$ is a fixed quantity throughout the manuscript. In particular, constants can possibly depend on the value of $p$. In this vein, we let
\begin{align} \label{theta}
  \th := \P(o \in \ccf)
\end{align}
denote the percolation probability, where $o = (0, 0)$ denotes the origin in $\Z^2$.

For a constant $\kappa \in (0,\theta)$ fixed throughout the paper, let  $\mc E^{\ms{uniq}}_n$ be the event that there exists a unique open connected ``giant'' component in $B(n) := [-n, n]^2 \cap \Z^2$, called as $\cc^n$, whose size is at least  $\kappa (2n)^2$  and all other components in $B(n)$ are of size at most $(\log n)^5$.
Then, by \cite[Theorem 4]{pp},   $\P(\mc E^{\ms{uniq}}_n) \rightarrow 1$ as $n\rightarrow \infty$.

~

In this work, we study supercritical percolation on $\Z^2$ from the angle of the isoperimetric profile. To formulate this precisely, we first say that a path $\g = (x_0,\cdots,x_n)$ in $ \Z^2$   (i.e.~all $x_i$s are lattice points) is a \emph{circuit} if  $x_n = x_0$, $|x_{k+1} - x_k |=1$ ($k=0,\cdots,n-1$) and it is simple in the sense that $x_0,\cdots,x_{n-1}$ are all distinct. 
Henceforth for a circuit $\g$ in $\Z^2$, $\vol(\g)$ denotes the collection of lattice points that are contained in $\g$ or inside the domain enclosed by $\g$. 
In addition, a path $\gamma$ in $\Z^2$ is called {open} (resp. {closed}) if it only consists of open (resp. closed) edges.

Under the event $\mc E^{\ms{uniq}}_n$, {denoting by $\mathcal C_n$ the family of open circuits in the giant component $\ccn$,} define the \emph{isoperimetric constant}
\begin{align} 	\label{eq:isop}
	\wtpn:= \inf \Big\{\f{|\g|}{|\ccn \cap \vol(\g) |}\co  \g \in \mathcal C_n,  
 \  |\vol(\g)|\le \frac12 |B(n)|\Big\},
\end{align}
where $|\cdot|$ denotes the cardinality of a set.
  We note that a slightly different definition of isoperimetry called \emph{Cheeger constant}, related  to the notion of edge-expander\footnote{A graph $G$ is called $\alpha$-\emph{edge-expander} if    $\frac{|\partial_G S| }{|S|}  \ge \alpha $ for all $S\su V(G)$ such that $\sum_{v\in S} d_G(v) \le 
 \frac{1}{2} \sum_{v\in V(G)} d_G(v)=|E(G)|$, where $\partial_G S$ denotes  the set of edges having exactly one endpoint in $S$.}, is used in \cite{blpr}, namely  
\begin{align}
\label{eq:blpr}
  \wt \Phi_n := \inf \Big\{ \frac{|\partial_{\ccf}U|}{|U|} : U \su \ccn,  \ 0< |U| \le \frac12 |\ccn|\Big\},
\end{align}
{where $\partial_{\ccf}U$ denotes the set of edges in $\ccf$ with exactly one endpoint in $U$ (note that $\ccn \su \ccf$ almost surely)}.  
Then, \cite[Theorem 1.4]{blpr} shows that there exists  a specific constant $\xi' > 0$ such that conditionally on  the event $\mc E^{\ms{uniq}}_n,$
\begin{align}
	\label{eq:wl}
	 \lim_{n\rightarrow \infty} n\wt \Phi_n  = \frac{\xi'}{\sqrt 2\th}   
\end{align}
in probability (see also \cite{gold2} for the LLN result
for the Cheeger constant ${|\partial_{\ccn}U|}/{|U|}$, i.e. the edge-boundary is taken in  $\ccn$  instead of  $\ccf$).

~

The  definition \eqref{eq:blpr} is natural from the graph-theoretic perspective, whereas our definition  \eqref{eq:isop} is inspired by the  geometric perspective:
The isoperimetric constant of  an  $n$-dimensional closed Riemannian manifold $M$  is defined by
\begin{align*}
    \inf_E  \frac{\text{Vol}_{n-1}(E)}{\min \{\text{Vol}_n(A),\text{Vol}_n(B)\}}
\end{align*} 
($\text{Vol}_k(\cdot)$  denotes the  $k$-dimensional volume measure),
where the infimum is taken over all $(n-1)$-dimensional submanifolds $E$ of $M$ which divide it into two disjoint submanifolds $A$ and $B$. Indeed in the definition \eqref{eq:isop}, instead of the edge-boundary $\partial_{\ccf}U$, the length of an open circuit $\g$ (called \emph{chemical distance}) is considered. 

The Cheeger constant \eqref{eq:blpr} has an intrinsic  nature and is thus related to the behavior (e.g. mixing time) of the random walk on the cluster. On the other hand  our notion of isoperimetry \eqref{eq:isop} is ``extrinsic'' and concerns the isoperimetric behavior of the cluster $\ccn$ with respect to the underlying space $\Z^2$. It incorporates the concept of chemical distance and density via the identity
\begin{align*}
         \f{|\g|}{|\ccn \cap \vol(\g) |} =\f{|\g|}{|\vol(\g)|}  \cdot \Big( \frac{|\ccn \cap \vol(\g)|}{|  \vol(\g) |} \Big)^{-1}
\end{align*}
(the second quantity above, without the inverse, is called \emph{density} of the cluster $\ccn$).  
{The chemical distance and cluster density are central observables in  the context of percolation, which have been extensively studied so far.} Indeed, their large deviation properties have attracted immense interest, due to the connection with Wulff construction (see Section \ref{sec 2.1} for details). As we will see  later, {the large deviation properties of these observables are incorporated into the variational problem \eqref{eq:isop} in a delicate way, which yields an interesting phase transition phenomenon for the isoperimetric constant.}

Note that the definitions \eqref{eq:isop} and \eqref{eq:blpr} are related to each other, as there is a  correspondence (albeit not exact) between an open circuit $\g \in \mathcal{C}_n$ and a region in $\ccn$. Indeed, for  $\g \in \mathcal{C}_n$, the associated region is $U=\ccn \cap \vol(\g)$. 
Conversely, for $U\su \ccn$, there exists $\g \in \mathcal{C}_n$ which traverses each edge at most once in each direction such that $U \su   \ccn \cap \vol(\g)$ and $| \{e\in \partial^+ \g : e \text{ is open}\}| \le |\partial_{\ccf}U|$, where $\partial^+ \g$ denotes the  ``right-boundary'' of $\g$ (see \cite[Eq. (2.1), (5.10)]{blpr} for explanations). {Two isoperimetric constants $\Phi_n$ and  $ \wt \Phi_n$  are typically 
 of order $n^{-1}$ and thus  equivalent (up to a multiplicative constant) in a typical regime.} 
Note that in the case of higher-dimensional lattices, the definition \eqref{eq:blpr} still makes  sense,  whereas in the definition \eqref{eq:isop} circuits $\g$ have to be  replaced with hypersurfaces. 

~

Our main result, Theorem \ref{thm:main} below,
concerns the decay rate for the probability of large deviation events for the isoperimetric constant $\wtpn$.
\begin{maintheorem}
\label{thm:main}
 The following large deviation estimates hold, where $\xi > 4$ is a specific constant (depending only on $p \in (\frac{1}{2},1)$)
which will be explicitly defined  in  \eqref{iso percolation} below.

 1. For any $t <2\sqrt2 $,
\begin{align} \label{theorem1}
\P(n\wtpn \le t \mid \mc E^{\ms{uniq}}_n) = 0.
\end{align}

2. For any $t\in (2\sqrt2,\frac{2\sqrt2}{\theta})$,
\begin{align} \label{theorem2}
-\ff < \liminf_{n\rightarrow \infty }\f1{n^{2}}\log \P(n\wtpn \le t \mid \mc E^{\ms{uniq}}_n) \le \limsup_{n \rightarrow \infty  }\f1{n^{2}}\log \P(n\wtpn \le t \mid \mc E^{\ms{uniq}}_n ) <0.
\end{align}

3. For any $t\in (\frac{2\sqrt2}{\theta}, \frac{\xi}{\sqrt 2\th})$,
\begin{align} \label{theorem3}
-\ff < \liminf_{n\rightarrow \infty }\f1n\log \P(n\wtpn \le t \mid \mc E^{\ms{uniq}}_n ) \le \limsup_{n\rightarrow \infty}\f1n\log \P(n\wtpn \le t \mid \mc E^{\ms{uniq}}_n ) <0.
\end{align}

4. For any $t>\frac{\xi}{\sqrt 2\th}$,
\begin{align} \label{theorem4}
-\ff < \liminf_{n\rightarrow \infty}\f1n\log \P(n\wtpn \ge t \mid \mc E^{\ms{uniq}}_n ) \le \limsup_{n\rightarrow \infty}\f1n\log \P(n\wtpn \ge t \mid \mc E^{\ms{uniq}}_n) <0.
\end{align}
\end{maintheorem}
In particular, Theorem \ref{thm:main} implies that   conditionally on  the event $\mc E^{\ms{uniq}}_n,$
\begin{align} \label{lln}
    \lim_{n\rightarrow \infty}n\wtpn = \frac{\xi}{\sqrt 2\th}
\end{align}
in probability. Note that this is an analog of the LLN result  for the Cheeger constant \eqref{eq:wl}, obtained in \cite[Theorem 1.4]{blpr}, for our notion of  isoperimetry \eqref{eq:isop}. In the following remark, we briefly elaborate on the comparison between \eqref{lln} and \cite[Theorem 1.4]{blpr}.

\begin{remark} [Comparison between $\Phi_n$ and   $\wt \Phi_n $] 

\label{remark 2.1}
It is shown \cite[Theorem 1.4]{blpr} that  the LLN constant $\xi'$ in  \eqref{eq:wl} for the Cheeger constant $\wt \Phi_n $ is given in terms of the isoperimetric problem w.r.t. the metric related to the edge boundary (called \emph{right-boundary distance}, see \cite[Section 2]{blpr}). One of the major contributions in \cite{blpr} is a verification of LLN for the right-boundary distance.

    In contrast, the LLN constant $\xi$ in \eqref{lln}  for the isoperimetric constant $\Phi_n$  is obtained in terms of  the isoperimetric problem w.r.t. the chemical distance. The chemical distance was already known to satisfy LLN (called \emph{time constant} in the FPP literature). Hence, with the aid of some geometric measure theoretic tools developed in  \cite[Section 4]{blpr}, it is not hard to deduce the LLN result \eqref{lln} for our notion of isoperimetry \eqref{eq:isop}. Our major contribution is the large deviation result for the isoperimetric constant, Theorem \ref{thm:main}. See the discussion after Remark \ref{remark 2.3} for the difficulties in the problem of large deviations for the isoperimetric constant. 
    
   While the Cheeger constant $ \wt \Phi_n$   and the isoperimetric constant $\Phi_n$  are both of  the order $n^{-1}$ 
    in the typical regime, their large deviation behaviors are  expected to be different from each other. 
     Indeed, regarding the Cheeger constant $ \wt \Phi_n$, we believe that  there will be \emph{no} regime in the lower tail  exhibiting  a  volume order large deviation. {Roughly speaking, this is because one can close most of the edges crossing the boundary of a macroscopic open circuit in the dual lattice, with a probability cost of at most $e^{-\Theta(n)}$. However,  in the upper tail regime, we still expect that the large deviation probability for the Cheeger constant is of  surface order.
     }

\end{remark}

\begin{remark}[Upper vs lower tail large deviations]
	When considering large deviations in random  structures, it is frequently the case that the upper and lower tail probabilities exhibit different orders  (see for example \cite{log,note,harel,dembin,ghn} for various types of large deviations problems in the context of  random graphs). For instance, in the case of 2D FPP with \emph{bounded} single-edge passage times, the large deviation probabilities for the  lower and upper tail of passage times are  $e^{-\Theta(n)}$  and  $e^{-\Theta(n^2)}$, respectively \cite[Theorems 5.2, 5.9]{kestAsp}. In fact, an exceptionally long connection from $o = (0, 0)$ to $ne_1 = (n ,0)$ is only possible when a large proportion of the edges in an $\Theta(n)$ neighborhood around the segment joining $o$ and $ne_1$ have atypically large passage times, which accounts for the probability cost $e^{-\Theta(n^2)}$  (see \cite{bgs,cz} for details). 

	In \cite[Open problem (6)]{blpr}, the large deviation question is raised whether this principle also holds for the isoperimetric constant in supercritical percolation. Theorem \ref{thm:main} shows that  this is not the case, indeed the upper tail large deviation  probability for the isoperimetric constant is of order $e^{-\Theta(n)}$, showing a different decay from the passage time in the FPP problem. Loosely speaking, this discrepancy can be explained as follows: For the chemical distance in supercritical percolation, the passage time of closed edges can be regarded as infinity, which violates the boundedness condition on a passage time, and the obstacle creating an atypically long chemical distance can be constructed by {planting $\Theta(n)$ closed edges} \cite{gm} (see Section \ref{sec 2.1.2} for detailed explanations).
\end{remark}

\begin{remark}[\emph{Multiple} phase transition in lower tail] \label{remark 2.3}
	Recall that for the passage time in FPP, we find a decay of order $e^{-\Theta(n)}$ for the entire lower tail probabilities that are compatible with the support of single-edge passage times \cite[Theorem 5.2]{kestAsp}. On the other hand, regarding the isoperimetric constant in supercritical percolation, we observe multiple phase transitions  in the lower tail, see parts 1, 2, 3 of Theorem \ref{thm:main}.  In fact,  the ``much'' lower tail regime (i.e. far from the typical value) exhibits  volume order large deviations $e^{-\Theta(n^2)}$, whereas the ``less'' lower tail regime  (i.e. closer to the typical value) exhibits  surface order large deviations  $e^{-\Theta(n)}$.   We refer to Section \ref{idea} for the detailed reasons of the appearance of multiple phase transitions.
\end{remark}

Note that the isoperimetric constant  is written as
\begin{align} \label{inc}
     \f{|\g|}{|\ccn \cap \vol(\g) |} =\f{|\g|}{|\vol(\g)|}  \cdot \Big( \frac{|\ccn \cap \vol(\g)|}{|  \vol(\g) |} \Big)^{-1}.
\end{align}
This shows that the isoperimetric constant   is closely related to the chemical distance and cluster density. Although we make use of the large deviation results for these observables \cite{ds,gm}, these results do not immediately yield large deviation properties on the isoperimetric constant, as isoperimetry is a variational problem determined by an intricate interplay between the density and chemical distances.  In fact, unlike other observables, the isoperimetric constant has a lower tail which exhibits a  phase transition, possessing both volume and surface order large deviations. Roughly speaking, as will be shown later, the ``much'' lower tail is only governed by the large deviation behavior of the cluster density, whereas the ``less'' lower tail is  mainly determined by the large deviation behavior of the chemical distance of an open circuit. {In this paper, we find the precise location in the lower tail where this phase transition phenomenon occurs.}

In order to analyze the large deviation behavior of the isoperimetric constant  in a precise way, we develop much more refined large deviation results controlling the cluster density in a uniform fashion and encoding geometric properties of geodesics in a quantitative sense. {This together with deterministic geometric measure theoretic arguments in \cite[Section 4]{blpr} will conclude the proof of Theorem \ref{thm:main}.} We refer   to Section  \ref{idea}  for elaboration on the proof strategy.

~

 In the following section, we briefly review the existing large deviation results regarding the cluster density and chemical distance in supercritical percolation.
 
  \subsection{Related works} \label{sec 2.1}
    There are various crucial observables in the context of supercritical percolation on  $\Z^d$, where the associated large deviation questions have attracted immense interest. In this section, we briefly elaborate on the large deviation principle regarding the density of clusters and the chemical distance in percolation, as they are closely related to the  isoperimetric constant in the perspective of \eqref{inc}.

\subsubsection{Density of clusters}
Density of the cluster $\cc$ in  $B(n)$ is defined to be $ \frac{|  \cc\cap B(n)|}{|B(n)|}$.
Recalling that $\theta$ denotes the probability that the origin is contained in the infinite  cluster $\ccf$,  by the law of large numbers and translation invariance, the density of $\ccf$ in a large square is close to $\theta$. Regarding  large deviations, the upper and lower tail of the density of $\ccf$ are shown to exhibit different orders of probabilities. Indeed, the upper tail probability has a volume order, i.e., $e^{-\Theta(n^d)}$, whereas the lower tail probability  has a {surface} order, i.e., $e^{-\Theta(n^{d-1})}$ (see \cite{ds} for $d=2$ and  \cite{gan} for $d\ge 3$).
 The surface order  lower tail large deviation is due to the existence of a hypersurface surrounding the origin, consisting only of   closed edges, which makes the cluster density atypically small. 
 
 The large deviation rate function for the cluster density has also been examined. Using the abstract large deviation framework \cite{leb} for the translation invariant FKG (Fortuin-Kasteleyn-Ginibre) measure on the particle systems,  \cite[Theorem 2]{ds} deduced the existence of   large deviation rate function   in the upper tail regime. Whereas in the lower tail regime, this framework is not applicable and thus identifying the  rate function becomes more delicate. In the seminal works  \cite[Theorem 6.1]{acc} ($d=2$) and  \cite[Theorem 2.12]{cerf1} ($d =3$), the  large deviation rate function  was shown to exist in the lower tail, which is characterized in terms of the surface tension.

  The surface order large deviation result, along with the existence of the  rate function, is essentially used  in the rigorous construction of the Wulff shape. For instance, using the large deviation result \cite[Theorem 6.1]{acc} for the cluster density in 2D supercritical percolation, it is shown \cite[Theorem 3]{acc} that 
conditionally on the event that the density of $\ccf$ is less than its typical value $\theta$, there is a large dual contour which is macroscopically close to the Wulff shape. 

We finally remark that there are analogous phase separation results for the low-temperature models in statistical mechanics such as Ising model. Indeed,
clusters in supercritical percolation, conditionally on having an atypically low density of $\ccf$,  correspond to the microcanonical ensemble, given an atypically small magnetization.
To be precise, consider the low-temperature Ising model  on $\Z^d$ with a plus boundary condition, conditioned that the empirical magnetization is less than the spontaneous magnetization. Then, there is a single droplet consisting of  minus spins surrounded by plus spins, whose macroscopic behavior is given by the Wulff shape.
The surface order large deviation for the empirical magnetization is   crucial  to establish this principle. See   \cite{ioffe1,ioffe2}  ($d=2$) and \cite{bod,cerf} ($d\ge 3$) for details.
 
    \subsubsection{Chemical distance}  \label{sec 2.1.2}
 Chemical distance in percolation is defined to be the minimum number of open edges required to connect two sites (denoted by $D(x,y)$, see \eqref{chemical}  for the precise definition).  
Due to the averaging effect, the chemical distance satisfies the LLN and has a deterministic macroscopic behavior \cite[Theorem 3.2]{gm2}. The limiting constant in this LLN is called \emph{time constant}   (see Proposition \ref{prop: time constant}).

Large deviation properties  of the chemical distance  have also been studied.
    In \cite{gm}, it was shown that  large deviation probabilities for the chemical distance $D(o,x)$ exhibit an upper bound $e^{-\Theta(|x|_1)}$ in both the upper and lower tail regimes. Also, the bound $e^{-\Theta(|x|_1)}$ also serves as a lower bound in the upper tail regime. 
To see this, for $x =e_1= (1,0,\cdots,0)$, we require that all  edges along the segment from  $-\lfloor \e n \rfloor e_1$ to 0 are open and all edges sharing exactly one endpoint with this segment
are closed, except for the edge connecting $- (\lfloor \e n \rfloor   +1) e_1$ and $-\lfloor \e n \rfloor  e_1$.  This yields the desired lower bound,  as all paths from $o$ to $ne_1$ are forced to visit $-\lfloor \e n \rfloor  e_1$ and the probability cost is $e^{-\Theta(n)}$.   Note that this result holds in \emph{any} dimension $d$.

 The lower tail regime is more  complicated and the large deviation speed  depends on the geometry of endpoints \cite{gm}: Denoting by $\mu$ the time constant, for any $x\in \Z^d$, for small enough $\e  
 = \e(x)>0,$
    \begin{align*}
    \frac{1}{n} \log \P( D(o,nx) \le (1-\e) n \mu(x)) 
    \begin{cases}
         = -\infty \quad &\text{if }\mu(x) = |x|_1, \\
         \ge (\log p)|x|_1 \quad &\text{if }\mu(x) > |x|_1.
    \end{cases}
    \end{align*}
    In addition, by the subadditivity, the above quantity has a limit (including $-\infty$) as $n\rightarrow \infty$, yielding the large deviation rate function.

    However, in the upper tail regime, it is highly non-trivial to identify the large deviation rate function for the chemical distance, due to the absence of submultiplicity.
    Very recently, when    $d\ge 3$, the  rate function is shown to exist  in the upper tail regime \cite{dembin}. Indeed, it is shown that the upper tail event occurs by the existence of space-time cut-points that force geodesics to become tortuous. The reason of delicacy in the case $d=2$ is that not only the presence of space-time cut-points but also more global scenarios may affect the upper tail events (see \cite[Section 1.3]{dembin}  for details). 

We emphasize that the upper tail large deviation problems for the chemical distance in supercritical percolation and for the passage time in FPP (with \emph{bounded} passage times) have significantly  different natures. Unlike the former case where the upper tail probability is of order {$e^{-\Theta(n)}$ regardless of the underlying dimension $d$}, in the latter case an atypically large passage time in FPP forces a global change of random environments, which accounts for the 
 probability of order $e^{-\Theta(n^d)}$.

\subsubsection{Percolation on the     complete graph} 

In this section, we mention   isoperimetric properties of    supercritical percolation on the complete graph. The size of the largest component in the Erd\H{o}s--R\'enyi graph  $G(n,p)$  undergoes a phase transition at $p_c=\frac{1}{n}.$ In subcritical regime ($p=\frac{1-\e}{n}$), critical regime ($p=\frac{1}{n}$) and supercritical regime ($p=\frac{1+\e}{n}$), the size of the largest
component is logarithmic, $\Theta(n^{2/3})$ and $\Theta(n)$ respectively, with high probability.
{It is shown \cite{er1}  that with high probability,  the 
giant component in the supercritical Erd\H{o}s--R\'enyi graph is an $\alpha$-edge-expander ($\alpha>0$ is a constant) decorated by paths
and trees of controlled size.}
We also refer to  \cite{er4,er2} for the alternative approach, where the  expansion property is deduced by developing a particular model contiguous to the giant component in  supercritical $G(n,p).$
 Note that these results  do not provide an LLN type result for the isoperimetric profile.

\subsubsection{Percolation on the  hypercube} We conclude this section by mentioning the isoperimetric result of supercritical percolation on the hypercube $Q_d:=\{0,1\}^d$.    The critical percolation parameter on the hypercube is $p_c=\frac{1}{d}$ \cite{cube}, i.e. when $p=\frac{1-\e}{d} $, all components have a size of order $O(d)$ and when $p=\frac{1+\e}{d} $, there is a unique giant component  whose size is linear in $n=2^d$ with high probability. Recently, it is shown  \cite{hypercube0} that in the supercritical regime,  the giant component is a $\beta d^{-5}$-expander\footnote{A graph $G$ is called an $\alpha$-\emph{expander} if $\frac{|N_G(S)|}{|S|}  \ge  \alpha$ for every $S \su V(G)$ with  $|S| \le  |V (G)|/2$, where $N_G(S)$ denotes the set of vertices in $S^c$ with a neighbor in $S$.} ($\beta>0$ is a constant) with high probability. As in the case of supercritical Erd\H{o}s--R\'enyi graphs, the results  \cite{hypercube0} do not provide an LLN result.

\subsection{Further directions} In this section, we state some interesting and important further questions related to our results.
\subsubsection{Existence of the rate function}
It is an interesting further direction to verify the  existence of the large deviation rate function for the isoperimetric profile.
% and then examine the conditional structure given the occurrence of  rare events.
To elaborate on this problem, let us consider the large deviation problem in FPP (with bounded passage times). Using a
subadditive argument,  Kesten \cite{kestAsp}
deduced the existence of the rate function in the lower tail regime. The existence of the rate function in the upper tail regime has remained open for a long time, until it is resolved in \cite{bgs} recently. Indeed, \cite{bgs} justified the presence of a certain inherent metric structure which accounts for the upper tail event.
{Note that although the result there is only stated in the case $d=2$, the analogous result holds in any dimension as well (see \cite[Section 1.1]{bgs} for the discussion).}
Recall also that for the chemical distance in supercritical percolation on $\Z^d$ with $d\ge 3,$ the  upper tail large deviation rate function is shown to exist  in the very recent work \cite{dembin}.

In the context of the isoperimetric profile in percolation, the large deviation probabilities do not enjoy submultiplicity, both in the lower and upper tail regimes. 
{We believe that the robust large deviation result for the chemical distance in 2D percolation is a crucial ingredient to prove the existence of the  rate function, which is currently out of reach.}

\subsubsection{Higher-dimensional lattice} \label{high}
It is also a natural  further  question to study large deviations for the isoperimetric constant in higher-dimensional lattices.
In \cite{gold}, it is shown that
the isoperimetric profile of the giant cluster in supercritical percolation on higher-dimensional lattices  $\Z^d$ ($d\ge 3$) satisfies LLN. Unlike $d=2$, in the higher-dimensional case $d\ge 3$, the macroscopic quantity encoding the size of boundary is defined as the infimum of {certain weights} of hypersurfaces (of dimension $d-1$), \emph{not} in terms of circuits, which is related to the minimal 
``cutsets'' in the FPP literature.

Recently, \cite{cut} developed the surface order large deviation theory for the cutsets in FPP for any dimension $d\ge 2$. This result, along with the arguments 
 in \cite{gold}, could be used to examine large deviations for the isoperimetric profile in high-dimensional percolation. We expect that similarly as in Theorem \ref{thm:main}, the upper tail exhibits a surface order large deviation, whereas there are multiple regimes in the lower tail exhibiting the volume and surface order large deviation probabilities respectively.

\subsection{Idea of proof} \label{idea}
In this section, we briefly provide key strategies to prove Theorem \ref{thm:main}.

\subsubsection{Lower bound on the upper tail}
We aim to construct a barrier consisting of $\Theta(n)$ closed edges, which makes the isoperimetric constant of any  open circuit arbitrary large. In particular, such barrier is constructed so that the enclosed region of \emph{any} circuit (not necessarily open), not intersecting with this barrier, is partitioned into several ``thin'' regions sharing only a small portion of boundary (see Figure \ref{fig} for the illustration). This makes every open circuit become tortuous and thus have a large isoperimetric profile. This argument can be made rigorous 
using the isoperimetric inequality in $\Z^2$ (see Lemma \ref{discrete iso}), applied to each partitioned  region.
 Since the probability cost of planting  $\Theta(n)$  closed edges is $e^{-\Theta(n)}$, this yields the desired lower bound on the upper tail of the isoperimetric constant.

\subsubsection{Upper bound on the upper tail}   
We show that with probability at least $1 - e^{-\Theta(n)}$, there exists an open circuit in $\ccn$ whose isoperimetric constant is not too large. Note that the isoperimetric constant can be written as
\begin{align} \label{1}
    \f{|\g|}{|\ccn \cap \vol(\g) |} =\f{|\g|}{|\vol(\g)|}  \cdot \Big( \frac{|\ccn \cap \vol(\g)|}{|  \vol(\g) |} \Big)^{-1}.
\end{align}
Hence, we aim to  upper bound ${|\g|}/{|\vol(\g)|}$ and lower bound the density of $\ccn$ in $\vol(\g)$. First, following the idea in  \cite[Proposition 5.5]{blpr},   the quantity ${|\g|}/{|\vol(\g)|}$  can be made as small as possible once the open circuit $\g$ is taken to be close to the boundary of  Wulff shape w.r.t. chemical distance  (see  Section \ref{sec:isop} for the definition).  One of our main contributions is a sharp  quantitative estimate on the probability of availability of such approximation. To accomplish this, one needs to understand the behavior of geodesics in a quantitative sense. Next, in order to control the second term in \eqref{1}, one has to control the density of $\ccn$ uniformly in all reasonable regions.

In fact, we establish the following principles:
\begin{enumerate}
  \item ``Uniform'' concentration phenomenon for the density of $\ccn$ (Proposition \ref{lem:cdens}). 
    \item ``Almost'' straight property of the ``almost''-optimal path  in $B(n)$ with respect to the chemical distance (Proposition \ref{cor:eop}).
\end{enumerate}
Indeed, we lower bound the  probability  of  each scenario above, which is sufficient to establish surface order large deviations for the upper tail of the isoperimetric constant.
Note that the upper and lower tails for the cluster density, in  the principle (1), exhibit different orders of large deviation probabilities, i.e. volume and surface order respectively.

\medskip

We briefly elaborate on each principle. Recall that the upper and lower tail large deviation probabilities for the density of $\ccf$ in $B(n)$ are of order $e^{-\Theta(n^2)}$ and $e^{-\Theta(n)}$ respectively \cite[Theorem 2]{ds}. As the isoperimetric constant is given by a variational problem, one has to control the density  of the cluster \emph{uniformly} in all admissible regions. In \cite[Lemma 5.3]{blpr}, a much weaker upper bound $e^{-c(\log n)^2}$ for the cluster density is obtained (both in the lower and upper tail regimes), uniformly over regions 
in $B(n)$. We  significantly improve this by establishing the correct ``uniform'' large deviation bounds. {This is done by partitioning the region  by small macroscopic boxes and then analyzing the cluster density in each box, depending on whether a box lies in the interior or boundary of the region.}

 Next, we establish a quantitative stability property of the geodesics in supercritical percolation, i.e. the principle (2). Using an idea in \cite[Proposition 3.2]{blpr}, we partition the segment connecting two points of macroscopic distance  and then take geodesics joining adjacent points. As the length of each geodesic is not large,
 the concatenation of these geodesics is (macroscopically) close to the line segment, provided that the number of partition is sufficiently large.  One can deduce that this occurs with probability at least $1 - e^{-\Theta(n)}$, using large deviation properties of the chemical distance.

\medskip

 Given the above principles, we conclude the proof as follows.   Consider the Wulff shape, a minimizer of the isoperimetry problem w.r.t. chemical distance. As an application of the  principle (2), we approximate  the (macroscopically rescaled) boundary of  Wulff shape by some open circuit $\g$. Using the  principle (1) (in particular the lower tail regime), one can show that the isoperimetric constant of $\g$ is ``almost'' at most   $\frac{\xi}{\sqrt{2}\theta}$. Since the exceptional probability of each principle above is of surface order, we establish the upper bound on the upper tail events.

\subsubsection{Lower bound on the lower tail}
First, the lower bound in \eqref{theorem2} is a straightforward consequence of the event of volume order  probability that all edges in $B(n)$ are open. Next, the lower tail in \eqref{theorem3} is implied by the event, of the surface order  probability, that the perimeter of a certain open circuit is atypically small. In fact, assume that all edges connected to $\partial B(\lfloor n/\sqrt{2} \rfloor)$ are open. Since  the density of the giant cluster in $B(\lfloor n/\sqrt{2} \rfloor)$  is asymptotically $\theta$, we deduce the desired lower bound by taking a circuit $\partial B(\lfloor n/\sqrt{2}\rfloor)$ in \eqref{eq:isop}.

\subsubsection{Upper bound on the lower tail}
First, one can deduce from \eqref{1} that regardless of a circuit, smallness of the isoperimetric constant in the ``much'' lower tail regime \eqref{theorem2} implies that the cluster density in the region enclosed by a circuit is atypically large. By the uniform upper tail large deviation result for the cluster density, which is of volume order (see  \eqref{upper} in  Proposition \ref{lem:cdens}), we obtain the upper bound in  \eqref{theorem2}.

In the ``less'' lower tail regime
\eqref{theorem3}, as in  \cite[Proposition 4.1]{blpr}, we approximate an open circuit $\g$   by a Jordan curve $\la$  in $\R^2$ such that $|\g| \ge (1-o(1)) \len_\mu(\la)$.  Our contribution is a quantitative estimate on the availability of such approximation, in fact we show that the exceptional probability is of surface order  (see Lemma \ref{lem:pr41}).  As $$\frac{|\g|}{|\vol(\g)|} \ge (1-o(1))\frac{|\len_\mu(\la)|} {| \intt (\la)|}$$ 
and the latter quantity can be lower bounded using the isoperimetric inequality  w.r.t. chemical distance, by applying the uniform concentration result for the  cluster density, we deduce the upper bound in  \eqref{theorem3}.

\subsection{Notations}
 For any subgraph $G$ of the integer lattice $\Z^2$, $V(G)$ and $E(G)$ denote the vertex and edge sets respectively. 
Throughout the paper, rectangles (and squares) are always assumed to have vertices only in $\Z^2$ and edges are parallel to the $x$ or $y$ axis. As a subgraph of the lattice $\Z^2$, the rectangle $R$ having corners $(n_1,m_1),(n_1,m_2),(n_2,m_1),(n_2,m_2)$ ($n_1\le n_2, m_1\le m_2$) has a vertex set $V(R) = [n_1,n_2]_\Z \times [m_1,m_2]_\Z$  (we denote $I_\Z:=I \cap \Z$ for any $I \su \R$) and an edge set $E(R) = \{ e\in E(\Z^2):\text{both of the  endpoints of $e$ belong to $V(R)$} \}.$
The size of this rectangle $R$  is defined to be $(n_2-n_1) \times (m_2-m_1).$

In addition, for any $x\in \R^2$ and $p\in [1,\infty) \cup \{\infty\}$,  $|x|_p$  denotes the $\ell^p$-norm of $x$.
For a subset $S \su \Z^2$, we let  $\text{diam}(S):=\max \{|x-y|_1: x,y\in S \}$ denote the diameter of $S$. We also follow the convention that the infimum of the empty set is $+\ff$. 
 
 {Finally, for non-negative functions $f(n)$ and $g(n)$, we use  the following Landau notations. We write $f(n) = O(g(n))$ if $\limsup_{n \rightarrow  \infty}f(n)/g(n) < \ff$ and  $f(n) = o(g(n))$ if $\lim_{n\rightarrow \infty}f(n)/g(n) = 0$. Furthermore, $f(n) = \Omega(g(n))$ means $g(n) = O(f(n))$ and $f(n) = \Theta(g(n))$ means $f(n) =  O(g(n))$ as well as $g(n)= O(f(n))$.}

\subsection{Acknowledgement}
The authors thank Barbara Dembin for useful observations on the various notions of isoperimetry. This helped us to improve the presentation of our paper.
KN's research is supported by  the National Research Foundation of Korea (NRF-2019R1A5A1028324, NRF-2019R1A6A1A10073887).

\section{Properties of 2D supercritical percolation} 
\label{sec:pre}
In this section, we state some useful properties of (bond) percolation in $\Z^2$ with a fixed supercritical connectivity probability $p>p_c=1/2$. We denote a \emph{cluster} to be a connected subgraph of the integer lattice $\Z^2$ consisting only of open edges. In particular, for a rectangle $R$,  we say that $\cc$ is a cluster in $R$ if it is connected and consists of open edges contained in $R$.

\subsection{Density of clusters}
In this section, we discuss a crossing property and the density of percolation clusters.
A cluster $\cc$ in a rectangle $  R$ is called \emph{crossing} if there exist open paths $\g_1$ and $ \g_2$ in $\cc$ connecting two parallel horizontal and vertical boundaries of $R$, respectively. In addition, the \emph{density} of the cluster $\cc$ in $R$ is defined to be $|\cc|/|R|$. 

It is known \cite{pp} that with high probability, there exists a unique crossing cluster in the square whose density is close to $\theta$ (recall that $\theta$ denotes the probability that the infinite cluster $\ccf$ contains the origin,  see \eqref{theta}), and the diameter of any other cluster is relatively small. From now on, we denote by $\sigma(E)$ the $\s$-algebra generated by any collection of edges $E \subseteq E(\Z^2)$.

\bet[Theorems 4 and 5 in \cite{pp} with $\phi_n := (\log n)^2$] \label{cluster in square}
For $\de>0$, 
let $\mc E_\de (B(n)) $ (measurable w.r.t.~$\sigma (B(n))$) be the event such that 
\begin{enumerate}
  \item There exists  a unique  cluster $\cc\su B(n)$ of diameter at least $(\log n)^2$.
  \item $\cc$ is crossing and satisfies
  \begin{align*}
	\frac{|\cc|}{|B(n)|} \ge \theta - \de .
  \end{align*}
\end{enumerate}
Then, there exists a constant $c>0$ such that for sufficiently large $n$,
\begin{align*}
  \P (\mc E_\de (B(n)))  \ge 1-e^{-c(\log n)^2}.
\end{align*}
\ent

Note that the condition (1) above in particular implies that the size of any other cluster than $\cc$ in $B(n)$ is at most $ 4(\log n)^4$.

For our further purpose, we need a generalization of  Theorem \ref{cluster in square} for supercritical percolation on general subgraphs $R_n$ of the lattice $\Z^2.$ For  given constants $a>b>0$, a series of rectangles $\{L_n\}_n$ in $\Z^2$ is called of \emph{reasonable} size (of scales $a,b$) if the side lengths of $L_n$ are less than $an$ and greater than $bn$. We consider a series of subgraphs $\{R_n\}_n$  of $\Z^2$ such that
\begin{enumerate}
    \item There are (possibly overlapping) rectangles $R_n^{\ell}$ ($\ell=1,2,\cdots,r_n$ with $r_n$ bounded in $n$) of reasonable size (of scales $a,b$)  such that $V(R_n) = \cup_\ell V(R_n^{\ell})$ and $E(R_n) = \cup_\ell E(R_n^{\ell})$.
    \item Consider a graph $G_n = (V(G_n),E(G_n))$ with $V(G_n)  = \{1,2,\cdots,r_n\}$ such that vertices $i$ and $j$ are connected by an edge if  and only if  $R_n^i\cap R_n^j$ is a rectangle of reasonable size (of scales $a,b$). Then, $G_n$ is connected.
\end{enumerate}
We define a slightly stronger  notion of crossing as follows. 
We say that a cluster $\cc$ in $R_n$ is \emph{strongly crossing} if  for every $\ell=1,2,\cdots,r_n$, for any interval $I_\ell$ with $|I_\ell| \ge (\log n)^3$ which is a subset of the boundary of $R_n^{\ell}$, $\cc$ intersects $I_\ell$.
\bepr[Existence and uniqueness of strongly crossing clusters] \label{cluster in rec}
Suppose that $\{R_n\}_n$ is a sequence of subgraphs of $\Z^2$ satisfying the above conditions for some scales $a>b>0$.
For $\de>0$, 
let $\mc E_\de (R_n) $ (measurable w.r.t. $\sigma (R_n)$) be the event  such that 
\begin{enumerate}
  \item There exists  a unique   cluster $\cc \su R_n$ of diameter at least $(\log n)^2$. 
  \item $\cc$ is strongly  crossing and satisfies 
  \begin{align} \label{12}
  \frac{|\cc|}{|R_n|} \ge \theta -\de .
  \end{align}
\end{enumerate}
Then, there exists a constant $c>0$ such that for sufficiently large $n$,
\begin{align*}
  \P (\mc E_\de (R_n))  \ge 1-e^{-c(\log n)^2}.
\end{align*}
\enpr

Since Proposition \ref{cluster in rec} can be obtained using the argument in the proof of  Theorem \ref{cluster in square}, we 
 present the proof in the Appendix.

~

By the LLN and translation invariance, the density of the infinite cluster $\ccf$ in a large square is asymptotically $\theta$. 
The large deviation probability for the density of $\ccf$ 
was established in \cite[Theorems 2 and 3]{ds}:
for any small constant $\de>0$, 
\begin{align} \label{ld}
	\P ( |\ccf \cap B(n)| \ge (\theta+\de) |B(n)|) = e^{-\Theta(n^2)}
\end{align}
and
\begin{align} \label{ld2}
	\P ( |\ccf \cap B(n)| \le (\theta-\de) |B(n)|) = e^{-\Theta(n)}.
\end{align}
As mentioned in the introduction, it is important to note that  the lower and upper tails exhibit different large deviation speeds.

~

The crucial ingredient in our quantitative analysis on the isoperimetric constant  $\wtpn$ is a version of \eqref{ld} and  \eqref{ld2} for the cluster $\ccn$ (\emph{not} for $\ccf$) \emph{uniformly} in all reasonable regions. Under the event  $\mc E^{\ms{uniq}}_n$, for any subset $S \subseteq B(n)$, define the density of $\ccn$ in $S$ as
\begin{align*}
    \theta_n (S) := \frac{|  \ccn\cap S|}{|S|}.
\end{align*} 
Then, for $\de> 0$, define
 \begin{align} \label{con1}
     s_{n, \de}^+  := \max_{\substack{\g: \textup{ circuit in } B(n)\\ \diam(\g) \ge \de n \\ n|\g|/|\vol(\g)| \le  \de^{-1}}}(\th_n(\vol(\g)) - \th)
 \end{align}
and
\begin{align} \label{con2}
    s_{n, \de}^-  := \max_{\substack{\g:\textup{ circuit in } B((1-\de)n) \\ \diam(\g) \ge \de n \\ n |\g|/|\vol(\g)| \le  \de^{-1}}}(\th - \th_n(\vol(\g))).
\end{align}

In the next proposition, we establish a  universal bound for the lower and upper tails on the density of $\ccn$, uniformly in all admissible regions. 

\bepr[Uniform concentration of density of clusters]
\label{lem:cdens}
For any small enough $\de>0$, 
\begin{align} \label{upper}
\limsup_{n \rightarrow \infty} \f1{n^2} \log\P(s^+_{n,\de} \ge  \de\mid \mc E^{\ms{uniq}}_n ) < 0 
\end{align} 
and 
\begin{align} \label{lower}
\limsup_{n \rightarrow\infty} \f1{n} \log \P(s^-_{n,\de} \ge  \de\mid \mc E^{\ms{uniq}}_n ) <0.
\end{align} 
\enpr

A crucial ingredient in the proof of 
 \eqref{upper} is the following  lemma which  provides a general upper bound on the probability of the existence of a cluster having an atypically high density in the rectangle $R$. 

\bel \label{lemma 3.4}
Let $\de>0$. Then,
there exist constants $M,c>0$ such that for any  rectangle $R$ whose side lengths are at least $M$,
\begin{align} \label{340}
	\P\big(\max_{\textup{cluster }\cc\su\Z^2} |\cc \cap R| \ge (\theta+\de) |R| \big) \le e^{-c|R|}.
\end{align}
\enl

This lemma is obtained by the similar reasoning as in \eqref{ld} and thus will be proved in the Appendix. Assuming Lemma \ref{lemma 3.4}, we  establish Proposition \ref{lem:cdens}.

\bep[Proof of Proposition \ref{lem:cdens}]
Since $\P(\mc E^{\ms{uniq}}_n )\rightarrow 1$ as $n\rightarrow \infty$, it suffices to prove that  for some  $c>0,$
\begin{align} \label{upperr}
 \P( \{s^+_{n,\de} \ge  \de\}  \cap  \mc E^{\ms{uniq}}_n ) \le e^{-cn^2}
\end{align} 
and 
\begin{align} \label{lowerr}
 \P( \{s^-_{n,\de} \ge  \de\} \cap  \mc E^{\ms{uniq}}_n ) \le e^{-cn}.
\end{align} 

We assume that the event  $\mc E^{\ms{uniq}}_n$ occurs so that the cluster $\ccn$ is well-defined.
Setting  $r  :=\lfloor \de^2  n\rfloor$, for $u\in \Z^2$, let $D_u := (2r)u + [-r, r-1]_\Z^2$. Then, the squares $D_u$s define a partition of $\Z^2.$
For any circuit $\g$ in $B(n)$ such that $\diam(\g) \ge \de n$, define
\begin{align*}
    S_{\text{int}}(\g) := \{u\in \Z^2 \co D_u \su \vol(\g) \sm \g\},\quad  S_{\text{bdr}}(\g) := \big\{u\in \Z^2 \co D_u \cap \g \ne \es\big\}.
\end{align*}  
Note that  $S_{\text{int}}(\g)$ and $S_{\text{bdr}}(\g)$ are disjoint. Also,  
\begin{align} \label{geo3}
     (2r)^2|S_{\text{int}}(\g)|+(2r)^2|S_{\text{bdr}}(\g)| \ge  |\vol(\g)|
\end{align}
and
\begin{align} \label{geometry0}
    (2r)^2|S_{\text{int}}(\g)| \le  |\vol(\g)|.
\end{align} 
Since $\diam(\g) \ge \de n$ and $r=\lfloor \de^2  n\rfloor$,  for any small $\de>0$ and  $u \in \Z^2$, the circuit $\g$ {cannot be entirely contained in the dilated square $\bar D_u:= (2r)u + [-2r, 2r-1]_\Z^2$}. Thus, for any $u \in S_{\text{bdr}}(\g)$, $\g$ contains at least two distinct lattice points each of which contained in $D_u$ and $\bar D_u^c$ respectively, which implies $|\g \cap \bar D_u| \ge |\g \cap  ( \bar D_u \sm D_u) | \ge r$. In addition, any lattice point $x$ in $\g$ can be contained in at most 9 squares $\bar D_u$s. Hence, 
\begin{align} \label{geometry}
    |S_{\text{bdr}}(\g)|r \le \sum_{u \in S_{\text{bdr}}(\g)} \big|\g \cap  \bar D_u \big| \le  9 |\g|.
\end{align}
With these preparations, we now establish \eqref{upperr} and \eqref{lowerr}.

~

\textbf{Proof of \eqref{upperr}.} 
{By Lemma \ref{lemma 3.4} together with a union bound,}  with probability at least  $1-e^{-\Theta(n^2)},$
\begin{align} \label{360}
	\max_{D_u \su B(n)}    \th_n(D_u)  < \th + \de.
\end{align} 
Since any $u\in S_{\text{int}}(\g)$ satisfies  $D_u \su B(n)$, by  \eqref{geometry0} and \eqref{geometry}, under the  event \eqref{360},
\begin{align} \label{geo2}
    |\ccn \cap  \vol(\g)| \le (\theta+  \de) (2r)^2  |S_{\text{int}}(\g)|+ (2r)^2|S_{\text{bdr}}(\g)| \le (\theta+  \de)   |\vol(\g)| + 36 r |\g|  .
\end{align}
Recalling $r  =\lfloor \de^2  n\rfloor$ and the condition $\de n |\g| \le |\vol(\g)|$, this implies that 
$
    \th_n(\vol(\g)) \le \th+40\de.
    $
    Since $\de > 0$ was arbitrary, we deduce \eqref{upperr}.  

\medskip

\textbf{Proof of \eqref{lowerr}.} 
{We use the following fact, established in  \cite{bm} (see \cite[Lemma 5.2]{blpr} for details): There exists $c'>0$ such that for any $1\le n'\le n$,
\begin{align*}
   \P( \ccn \cap B(n') = \ccf \cap B(n') )  \ge 1-e^{-c'(n-n')}.
\end{align*}
In particular,  $\ccn \cap B((1 - \de/2)n) = \ccf \cap B((1 - \de/2)n)$ with probability at least $1-e^{-\Theta(n)}.$}
Under this event, for any  $D_u \su B((1 - \de/2)n),$ $$\th_n(D_u) =\frac{|\ccn \cap D_u|}{|D_u|} = \frac{|\ccf \cap D_u|}{|D_u|}. $$  Hence, \eqref{ld2} together with a union bound imply that with probability at least $1-e^{-\Theta(n)}$,
 \begin{align}   \label{361}
	\min_{D_u \su B((1 - \de/2)n)}    \th_n(D_u)  > \th - \de.
\end{align}
Since $\g$ is contained in $B((1-\de)n),$  for small $\de>0$, any $u\in S_{\text{int}}(\g)$ satisfies  $D_u \su B((1 - \de/2)n)$. 
Therefore, using the lower bound on $|S_{\text{int}}(\g)|$ obtained by \eqref{geo3} and  \eqref{geometry}, under the  event \eqref{361},
 \begin{align*}
     |\ccn \cap \vol(\g)| \ge (\theta-\de) (2r)^2  |S_{\text{int}}(\g)|  &\ge   (\theta-\de)  (2r)^2  \Big(\frac{|\vol(\g)|}{(2r)^2} -  \frac{9|\g|}{r}\Big) \\
     & \ge   (\theta-\de)   |\vol(\g)| \Big(1-\frac{36r}{\de n} \Big) \ge (\theta-40\de)   |\vol(\g)|,
\end{align*}
where  we used the condition $|\vol(\g)| \ge \de n|\g|$ in the second last inequality.   By the arbitrariness of $\de>0$, we deduce \eqref{lowerr}. 
\enp

\subsection{Chemical distance} \label{sec 3.2}
In this section, we introduce some useful properties of the chemical distance in supercritical percolation.  We denote by $|\g|_1$ the  length of a lattice path $\g$ in $\Z^2$.  Then, define the \emph{chemical distance} $D: \Z^2 \rightarrow \Z $ by
\begin{align} \label{chemical}
  D(x,y):= \inf_\g |\g|_1,
\end{align}
where the infimum is taken over  all open paths $\g$ joining two points  $x$ and $y$. We set $D(x,y):=\infty$ if $x$ and $y$ are not connected through open edges. Obviously, the chemical distance is always at least the $\ell^1$-distance.

Although the chemical distance $D(\cdot,\cdot)$ is a random distance, it has a deterministic macroscopic behavior due to the averaging effect. For $x\in \Z^2$, we denote by  $\{x \longleftrightarrow \infty\}$ the event that $x$ belongs to the infinite cluster $\ccf$.

\bepr[Theorem 3.2 in \cite{gm2}] \label{prop: time constant}
Consider supercritical    percolation in $\Z^2.$ Then,
there exists a norm $\mu$
on $\R^2$ such that almost surely on the event $\{o \longleftrightarrow \infty\}$,
\begin{align*}
  \lim_{n\rightarrow \infty} \frac{D(o,T_{n,x}x)}{T_{n,x}} = \mu(x),\quad \forall x\in \Z^2,
\end{align*}
where $(T_{n,x})_{n\ge 1}$ is any increasing sequence of positive integers $k$ such that $\{kx \longleftrightarrow \infty\}$.

\enpr
The norm $\mu$ above is called \emph{time constant}.
 Although this proposition was established  for supercritical percolation in general lattices $\Z^d$  \cite[Theorem 3.2]{gm2}, we only state the two-dimensional case for our purpose. Since the chemical distance is always at least the $\ell^1$-distance, we have 
\begin{align}
  \mu(x) \ge |x|_1,\quad \forall  x\in \R^2.
\end{align}
However, one can improve this trivial inequality and establish that the chemical distance uniformly (strictly) dominates the $\ell^1$-distance in the ``almost'' horizontal and vertical directions. This is a key reason behind a  {multiple} phase transition phenomenon for the lower tail large deviations of the isoperimetric profile.

%
%DOM PROP
%
\bepr[Domination for almost horizontal and vertical directions] \label{prop: time constant dominance}
There exist constants $\tau,\zeta>0 $  such that the following holds. For any $ x=(x_1,x_2)\in \Z^2$ such that $|x_2| \le \tau |x_1| $ or $ |x_2| \ge \frac1{\tau}|x_1| $,
\begin{align} \label{351}
  \mu( x) \ge (1+\zeta) |x|_1.
\end{align}
\enpr

The motivation behind the ``almost'' horizontal and vertical condition on $x\in \Z^2$ is the following. Under such condition, the number of paths in $\Z^2$ (not necessarily open) connecting $o$ and $x$, whose $\ell^1$-length is less than $(1+\de )|x|_1$, is at most $e^{c |x|_1}$ for some
 $c = c(\tau,\de)>0$. This constant $c>0$ can be made arbitrary small, provided that $\tau,\de>0$ are small enough. By a union bound,
\begin{align}
	\P \big(\min \{ |\g|_1\co \g \text{ open  path connecting $o$ and $x$}\}<(1+\de)|x|_1\big) \le e^{c|x|_1} p^{|x|_1}.
\end{align}
This bound converges to 0 as $|x|_1 \rightarrow  \infty$ if $c>0$ is small enough, i.e. when   $\tau,\de>0$ are small enough. Therefore, {by the definition of  the norm $\mu$,} we obtain \eqref{351}.
 The rigorous proof will be presented in the Appendix (see Section \ref{sec 7.3}).

~

Finally, we state a large deviation bound for the chemical distance, obtained in 
\cite[Theorem 1.2]{gm}. It claims that both lower and upper large deviation probabilities for  the chemical distance  are at most exponential.
\begin{proposition}[Theorem 1.2 in \cite{gm}] \label{ld distance}
For any constant $\de>0,$  
\begin{align} \label{dist}
   \limsup_{|x|_1 \rightarrow \infty} \frac{1}{|x|_1}\log  \P\Big(o \longleftrightarrow x, \Big\vert  \frac{D(o,x)}{\mu(x)}  - 1 \Big\vert  \ge  \de \Big) < 0. 
\end{align}
Here, $\{o \longleftrightarrow x\}$ denotes   the event  that $o$ and $x$ are connected through open edges.
\end{proposition}

To obtain a version of  \eqref{dist} for all $x\in \Z^2$, not only for vertices connected to the origin, we define  the notion of ``closest'' vertex in the cluster $\ccn$. For the moment,  we assume that the event $\mc E^{\ms{uniq}}_n$ occurs so that $\ccn$ is well-defined.   For $R>0$, define a continuous version of the square of size $2R$ as $B^{\text{con}}(R) := \{x\in \R^2 : |x|_\infty \le R\}.$ Let $\{\eta_z\}_{z\in \Z^2}$ be i.i.d.~uniform random variables on $[0,1]$ independent of everything else. 
Under the event $\mc E^{\ms{uniq}}_n,$ for a possibly non-lattice point
$x\in B^{\text{con}}(n)$, define $[x]$ to  be the vertex in $\ccn$ which is nearest to $x$ in the $\ell^\infty$-norm, and taking the one with the smallest $\eta_z$ if there is a tie  (although $[x]$  may depend on $n$, we omit the notation of $n$-dependency for the sake of  readability).

Then, we claim that $[x]$ is close to $x$ with high probability, i.e., there exists  $c>0$ such that for any $m> (\log n)^8$ and $r> 0$, 
\begin{align} \label{311}
	  \max_{x\in B^{\text{con}}(n-m)}   \P( | [x]-x|_\infty > r \mid \mc E^{\ms{uniq}}_n) \le e^{-c \min \{r, m-(\log n)^6\}}.
\end{align}
Note that a variant of the above estimate is proved in \cite[Lemma 2.7]{blpr}: denoting by  $[x]_\infty$ the closest vertex of $x$ in  $\ccf$,
$$  \P( | [x]_\infty-x|_\infty > r)  \le e^{-cr}.$$
A similar argument, which we present now, can be exploited to deduce \eqref{311}. 
 By monotonicity of the event in \eqref{311}, we may assume that $r \le m- (\log n)^6$.  We show that for $x\in B^{\text{con}}(n-m),$  the event $\{ |[x]-x|_\infty > r \} \cap \mc E^{\ms{uniq}}_n$ implies that there is no crossing between two boundaries of the annulus $B(n) \setminus \{ y\in \Z^2: |y-x|_\infty < r \}$. First, any cluster in $B(n)$, other than  $\ccn$, cannot  cross since this cluster has the size at most $(\log n)^5$ and the $\ell^\infty$-distance between two boundaries is at least $ m-r\ge (\log n)^6$. Also,  as $\{ |[x]-x|_\infty > r \}$, $\ccn$ cannot intersect the inner boundary of the annulus, $ \{ y\in \Z^2: |y-x|_\infty = r \}$.
 
 Hence, there exists an open circuit in a dual  percolation which is subcritical,  enclosing the $(2r\times 2r)$-square centered at $x$. Since the probability of such an event is exponentially decaying in $r$, we deduce \eqref{311}.

~

  From \eqref{311}, we deduce that for any    $\e,\iota>0$, there exists $c>0$ such that for any $x,y \in B^{\text{con}}(  (1 - \e)n  )$, conditionally on $\mc E^{\ms{uniq}}_n$, with probability at least {$1-e^{-cn}$,} 
\begin{align} \label{310}
   |([y] - [x]) - (y - x)|_1 /2 \le  |([y] - [x]) - (y - x)|_\infty \le |[y] - y |_\infty + |[x] - x |_\infty   \overset{\eqref{311}}{\le} \iota n.
\end{align}
From this along with    Proposition \ref{ld distance}, one can estimate $D([x],[y])$ as follows (note that the geodesic connecting $[x],[y]\in \ccn \su B(n)$ can leave $B(n)$): For any $\e,\e'>0$, there exists $c>0$ such that 
  for any $x,y \in B^{\text{con}}(  (1 - \e)n  )$ with  $|y-x|_1 \ge  \e' n,$ 
\begin{align} \label{31200}
  \P(
 | D([x],[y]) - \mu([y] - [x]) | \ge   \e  \mu([y] - [x])  \mid \mc E^{\ms{uniq}}_n) \le e^{-cn}.
\end{align}
To see this,   by \eqref{310} with small enough $\iota>0$, 
\begin{align*}
    \P(|[y] - [x]|_1 \ge  \e'n/2) \ge  1- e^{-cn}.
\end{align*}
Hence, \eqref{31200} follows by Proposition \ref{ld distance} together with {translation invariance and a union bound over all possible (random) points $[x]$ and $[y]$ in $B(n)$}.

%
%SEC ISOP
%
\section{Isoperimetric inequality}
\label{sec:isop}

In Sections \ref{ss:41} and \ref{ss:42}, we introduce some basic aspects of the isoperimetric inequality in the continuous and discrete spaces, respectively.  By a Jordan curve in $\R^2$, we mean a curve that is simple and closed.

%
%SS 41
%
\subsection{Isoperimetry in $\R^2$}
\label{ss:41}
We introduce the isoperimetric inequality w.r.t.~a general norm $\rho$ on $\R^2$.  
For any  curve $\la$ in $  \R^2$, i.e., a continuous map $\la: [0,1] \rightarrow \R^2$, define the \emph{$\rho$-length} of $\la$ to be
\begin{align*}
  \len_\rho (\la):= \sup_{N \ge 1}\sup_{0\le t_1 < \cdots < t_N \le 1} \sum_{i=1}^N \rho (\la(t_{i+1}) - \la(t_i) ) .
\end{align*}
The curve $\la$ is called \emph{rectifiable} if $\len_\rho (\la)<\infty$  for any norm $\rho$ on $\R^2$. For a Jordan curve $\la$ in $\R^2$,  we denote by $\intt(\la)$ the {interior of $\la$}, i.e., the unique bounded component of  $\R^2 \sm \la.$

The \emph{isoperimetric constant} w.r.t. the norm $\rho$ is defined to be  
\begin{align} \label{iso}
  \Iso_\rho:= \inf \{\len_\rho (\la): \la \ \text{is a Jordan curve in} \ \R^2, \
|\intt(\la)| = 1\},
\end{align}
where $|\cdot|$ denotes the Lebesgue measure on $\R^2.$

According to \cite{blpr,wulff}, the minimizer of \eqref{iso} can be characterized as follows. Define the region (called \emph
{Wulff shape})
\begin{align} \label{wulff}
 W_\rho:=  \cap_{\textbf{n}\in \R^2: |\textbf{n}|_2=1 } \{x \in \R^2 : \textbf{n} \cdot x \le \rho (\textbf{n})\}.
\end{align} 
Then, the boundary of $W_\rho$, normalized by its area, becomes the minimizer of the variational problem  \eqref{iso}. In addition, $W_\rho$ is a convex and compact set containing the origin. For the further purpose,
it is also useful to note that $ W_\rho$ is the unit ball w.r.t. the dual norm $\rho^*$  of $\rho,$ defined as
\begin{align}\label{dual}
     \rho^*(y):= \sup  \{x\cdot y: x\in \R^2, \rho(x) \leq 1 \}.
\end{align}
Furthermore, the minimizer of \eqref{iso} is unique up to translation \cite{taylor}.

~

In particular, in the special case when $\rho$ is the $\ell^1$-norm $| (x_1,x_2)|_1 := |x_1|+|x_2|$,  the isoperimetric constant is equal to 4, i.e., $  \Iso_{|\cdot|_1}=4$. 

\bel[Isoperimetric inequality w.r.t. $\ell^1$-norm] \label{lemma 41}
It holds that
\begin{align}\label{iso: l1 norm}
\inf_{\la} \frac{\len_{|\cdot|_1} (\la)}{\sqrt{|\intt(\la)|}} \ge 4,
\end{align}
where the infimum is taken over all Jordan curves $\la$ in $    \R^2$.
The equality is obtained if and only if $\lambda$ parametrizes a boundary of the square whose edges are parallel to the $x$ or $y$ axis.
\enl
\begin{proof}
By the characterization of the Wulff shape in \eqref{wulff}, 
$$W_{|\cdot|_1} = B^{\text{con}}(1) = \{(x_1,x_2) \in \R^2: -1\le x_1\le 1, -1\le x_2\le 1\}.$$ 
Let $\lambda_{|\cdot|_1} $ be a Jordan curve which parametrizes (say, clockwise) $ \partial B^{\text{con}}(1/2)$ in a constant speed, the boundary of the normalized Wulff shape $B^{\text{con}}(1/2)$. Since $\len_{|\cdot|_1} (\lambda_{|\cdot|_1})=4$, we deduce that $  \Iso_{|\cdot|_1}=4$. Therefore, by a scaling argument,
we obtain \eqref{iso: l1 norm}.
 The equality condition in \eqref{iso: l1 norm} follows from 
the uniqueness of the  minimizer of \eqref{iso} up to translation.
\end{proof}

Note that unlike $\len_{|\cdot|_2} (\la)$, the length $\len_{|\cdot|_1} (\la)$ is not invariant under the rotation. For instance, the rotated Jordan curve $\lambda_{|\cdot|_1}'$ which parametrizes $\{(x_1,x_2)\in \R^2: |x_1| + |x_2| = 1/\sqrt2 \} $  satisfies $\len_{|\cdot|_1} (\lambda_{|\cdot|_1}') =4\sqrt2\neq \len_{|\cdot|_1} (\lambda_{|\cdot|_1}) =4 $.

%
%SSSEC ISOP
%
\subsubsection{Isoperimetry w.r.t. chemical distance}

Now, we consider the isoperimetric constant w.r.t.~the norm $\mu$, a time constant in supercritical percolation.
Let
\begin{align} \label{iso percolation}
\xi:=  \Iso_\mu:= \inf \{\len_\mu (\la): \la \ \text{is a Jordan curve in $\R^2$ with} \
|\intt(\la)| = 1\}
\end{align}
(recall that  the norm $\mu$ is    defined in Proposition \ref{prop: time constant}).

In the next lemma, with the aid of Proposition \ref{prop: time constant dominance}, we show that the isoperimetric constant w.r.t. the chemical distance, $\Iso_\mu$, is strictly greater than the isoperimetric constant w.r.t. $\ell^1$-norm, $\Iso_{|\cdot|_1}$. This accounts for the \emph{multiple} phase transition in the lower tail large deviation  for the isoperimetric constant.

\bel[Strict domination in isoperimetric constant]
Let $\mu$ be the time constant defined in Proposition \ref{prop: time constant}. Then,
\begin{align*}
\textup{Iso}_\mu>\textup{Iso}_{|\cdot|_1}=4.
\end{align*}
\enl

\bep
Let $\lambda_{|\cdot|_1} $ be a Jordan curve parameterizing  $ \partial B^{\text{con}}(1/2)$  and
 $\la_\mu$ be a Jordan curve attaining the minimum of the variational problem \eqref{iso percolation}, unique up to translation. We consider two different cases, depending on whether $\la_\mu = \lambda_{|\cdot|_1}$ (whose traces in $\R^2$ coincide up to translation)  or not.

In the case $\la_\mu = \lambda_{|\cdot|_1}$, by Proposition \ref{prop: time constant dominance},
\begin{align*}
  \len_\mu (\la_\mu) =  \len_\mu (\la_{|\cdot|_1})  \ge (1+\zeta) \len_{|\cdot|_1}(\la_{|\cdot|_1})=4(1+\zeta),
\end{align*}
where in the above inequality we used that  $\lambda_{|\cdot|_1}$ consists  only of segments parallel to the $x$ or $y$ axis.

Whereas, in the case $\la_\mu \neq \lambda_{|\cdot|_1}$,
\begin{align*}
  \len_\mu (\la_\mu) \ge \len_{|\cdot|_1}(\la_\mu) >4,
\end{align*}
where the first inequality follows from the fact $\mu(x) \ge |x|_1$   and the last inequality follows from the uniqueness of the Wulff shape (w.r.t. $\ell^1$-norm) up to translation. 
\enp

%
%SS 42
%
\subsection{Isoperimetry in $\Z^2$}
\label{ss:42}

In this section, we state a discrete version of the isoperimetric inequality  on $\Z^2$. Recall that for a circuit $\g$ in $\Z^2$, $\vol(\g)$ denotes the collection of lattice points that are on   $\g$  or inside the domain enclosed by $\g$.

\bel[Isoperimetric inequality in $\Z^2$] \label{discrete iso}
For any $\e>0$, there exists $R>1$ such that the following holds. For any circuit  $\g$ in $\Z^2$ satisfying $|\vol(\g)| > R$ and $|\g| \le |\vol(\g)|^{2/3}$,
\begin{align} \label{371}
  |\g| \ge (4-\e) \sqrt{|\vol (\g)|}.
\end{align}
In addition, there exists a constant $c_0>0$ such that for \emph{all} circuits $\g$ in $\Z^2$,
\begin{align} \label{372}
  |\g| \ge c_0\sqrt{|\vol (\g)|}.
\end{align}
\enl

Note that for any  circuit $\g$ in $\Z^2$,  we have $|\g|=|\g|_1$ (i.e. the number of vertices is equal to the length of a circuit).
Before proving  Lemma \ref{discrete iso}, we present a useful inequality relating $|\vol(\g)|$ and $|\intt(\g)|  $ for a circuit $\g$ in $\Z^2$, where $\intt(\g)$ denotes the region enclosed by   $\g$ regarded as a Jordan curve in $\R^2$ by following the edges in $\g$ at linear speed.
Let $I$ and $B$ be the number of lattice points interior to $\g$ and on $\g$ respectively. Then, $|\g| = B$ and $    |\vol(\g)| = I+B$.
Also,
by Pick's theorem \cite{pick},
$
|\intt(\g)|  = I+ \frac{B}{2}-1.  
$
Hence, we deduce that
\begin{align} \label{diff}
\big||\vol(\g)|  - |\intt(\g)|  \big| \le  |\g|.
\end{align}

\begin{proof}[Proof of Lemma \ref{discrete iso}] 
We identify a circuit  $\g$ in $\Z^2$ with a Jordan curve   in $\R^2$.
Then, by a continuous version of the isoperimetric inequality (see Lemma \ref{lemma 41}), 
\begin{align*}
|\g| = \len_{|\cdot|_1} (\g) \ge 4 \sqrt{|\intt (\g)|} \overset{\eqref{diff}}{\ge} 4 \sqrt{ |\vol(\g)| - |\g|}   \ge (4-\e) \sqrt{|\vol (\g)|},
\end{align*}
where the last inequality follows from the conditions imposed on $\g$ with a large enough $R$. Hence,
we deduce \eqref{371}.

To prove \eqref{372}, note that  \eqref{371} with $\e=1$ (along with the corresponding  $R$) implies that, for any circuit $\g$ such that $|\vol(\g)| > R$ and $|\g| \le |\vol(\g)|^{2/3}$, we have $ |\g| \ge 3\sqrt{|\vol (\g)|}$. In addition, since there are only finitely many circuits $\g$ in $\Z^2$ (up to translation) with $|\vol(\g)| < R$,  there exists  a constant $c_1 = c_1(R)>0$ such that $  |\g| \ge c_1 \sqrt{|\vol (\g)|}$ for any such circuit $\g$. On the other hand, if the circuit $\g$   satisfies $|\g| > |\vol(\g)|^{2/3}$, then
\begin{align*}
  |\g| \ge |\vol(\g)|^{2/3} \ge \sqrt{|\vol (\g)|}.
\end{align*}
Taking $c_0:= \min \{c_1,1\}$, we establish \eqref{372}.

\end{proof}

The inequality \eqref{372} provides a straightforward lower bound on the isoperimetric constant. In fact, for any    circuit $\g$ in $\Z^2$ satisfying $|\vol(\g)| \le  \zeta n^2$ ($\zeta>0$ is a constant),
\begin{align} \label{333}
  n \frac{|\g|}{|\vol(\g)|} \ge  \frac{c_0 n}{\sqrt{|\vol(\g)|}} \ge \frac{c_0}{\sqrt{\zeta}}.
\end{align}

%
%SEC UP
%
\section{Upper tail large deviations}
\label{sec:up000}
\subsection{Lower bound on the upper tail}
\begin{figure}
	\centering
 \begin{tikzpicture}
	 \draw (-5,-5) rectangle (5,5);

	 \fill[red!20,opacity=.7] (-5,.5) rectangle (5,3);
	 \fill[green!20,opacity=.6] (3.6,-5) rectangle (5,5);
	 \fill[white,opacity=1] (0.4,1.8) rectangle (0.8,2.2);
	 \draw[green, thick, dashed] (3.6,-5) rectangle (5,5);
	 \draw[red, very thick, dash dot] (-5,.5) rectangle (5,3);

	  \foreach \x in {-50,...,35}{
	 \draw[black,  thick] (\x * 0.1, 3) -- (\x * 0.1,3.5);
	 \draw[black,  thick] (\x * 0.1, 0) -- (\x * 0.1,.5);
	 \draw[black,  thick] (\x * 0.1, -3) -- (\x * 0.1,-2.5);
 }

	   \draw[blue, very thick] plot[smooth, tension = 0.3] coordinates{(4.5,-4.7)(0, -4) (-4.2,-4.8) (-4.2,-3.1) (0, -3.5) (4.0,-3.3) (4.0,-2.0) 
	 (0, -1.5)  (-4.2,-2.0) (-4.2,-0.2) (0, -0.5)(4.2,-0.2) (4.5, 0.5) };
	 \draw[blue, very thick] (4.5,0.5)--(4.3,0.5);
	 \draw[blue, very thick] (4.1,0.5)--(3.9,0.5);
	   \draw[blue, very thick] plot[smooth, tension = 2.3] coordinates{ (4.3, 0.5) (4.2, 0.2) (4.1,0.5)};
	  \draw[blue, very thick] plot[smooth, tension = 0.3] coordinates{(3.9,0.5)(0, 1.3) (-4.2,0.8) (-4.2,2.8) (0, 2.5) (4.3,2.8) (4.4, 3.8) (4.1, 3.0)};
	 \draw[blue, very thick] (4.1,3.0)--(3.8,3.0);
	\draw[blue, very thick] plot[smooth, tension = 0.3] coordinates{ (3.8, 3.0) (3.8, 4.0) (0, 4.2) (-4.2,3.8) (-4.2,4.7) (0, 4.4)(4.5,4.7)  (4.5,-4.7)
	   };
	   \coordinate[label=-180:{\textcolor{red}{$R_i$}}] (A) at (0.9, 2);
	   \coordinate[label=-180:{\textcolor{green}{$R$}}] (A) at (4.6, -0.9);
	   \coordinate[label=-180:{\textcolor{blue}{$\g$}}] (A) at (0.0, -1);
	   \coordinate[label=-180:{{$\lfloor \t n\rfloor$}}] (A) at (4.8, -5.4);
	   \coordinate[label=-180:{{$\lfloor  n/k\rfloor$}}] (A) at (-5.2,1.8);
    \draw [decorate,    decoration = {brace,mirror}] (3.6,-5.1) --  (5,-5.1);
    \draw [decorate,    decoration = {brace}] (-5.1,.5) --  (-5.1,3);
 \end{tikzpicture}
 \caption{Circuit in $\Tilde{B}_n$}
 \label{fig}
\end{figure}
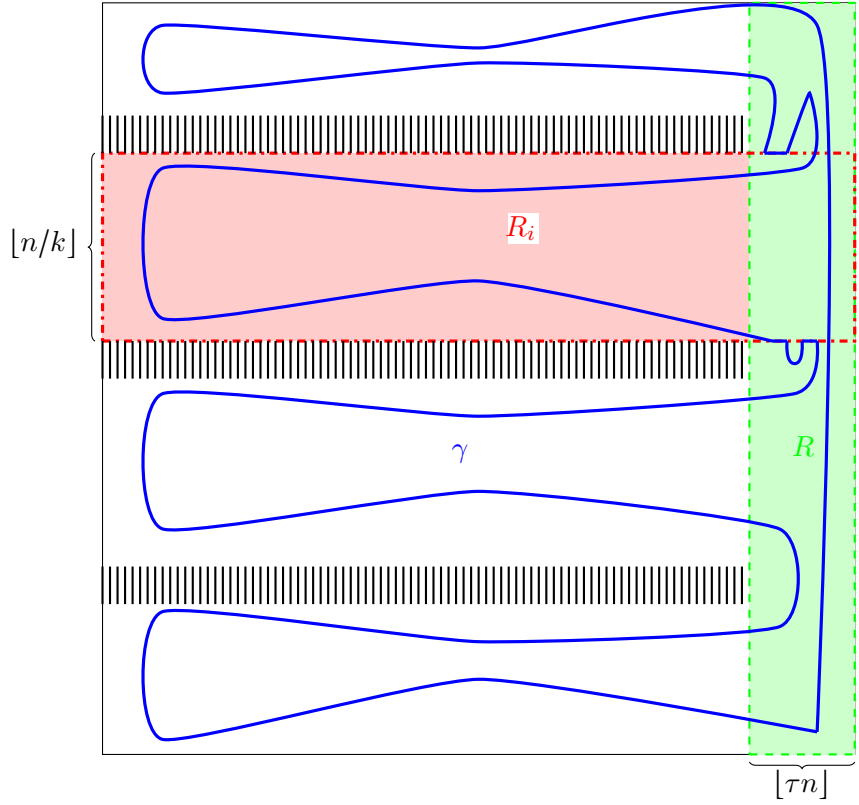

In this section, we establish the lower bound in \eqref{theorem4}, by verifying that for any constant $t>0$,
\begin{align}
  \liminf_{n \rightarrow \infty} \frac{1}{n}  \log\P (n\wtpn \ge t\mid \mc E_n^{\ms{uniq}}) > -\infty.
\end{align}
As mentioned in Section \ref{idea}, we construct a collection of barriers consisting of $\Theta(n)$ closed edges.
For $v\in \Z^2$, we denote by $v_1$ and $v_2$ the $x$ and $y$ coordinate of $v$, respectively.
Let $\tau>0$ be a small constant and $k\in \N$ be a large integer which will be chosen later. For each $i = -k+1,\cdots,k-1$, define the collection of edges
\begin{align}
  E_i:= \{e=(v,w) \in E(\Z^2) : v_1=w_1\in [-n,n- \lfloor \tau n \rfloor ]_\Z, \ v_2=w_2- 1=i \lfloor n/k \rfloor\},
\end{align}
and set $\bar{E}: = \cup_{i=-k+1}^{k-1} E_i$ (i.e. $\bar{E}$ is the set of vertical edges in Figure \ref{fig}).
Edges in $\bar{E}$ will play a role of a barrier, which enforces the isoperimetric constant of any circuit not intersecting with $\bar{E}$ to be  large. Define the event
\begin{align*}
  \cA_n:= \{ \text{All edges in} \ \bar{E}\ \text{are closed}\}
\end{align*}
(we suppress the notation of $k,\tau$-dependence  in $\cA_n$ for the sake of readability). 
Since $|\bar{E}| = \Theta(n)$,
\begin{align} \label{event a}
  \P(\cA_n ) = e^{-\Theta(n)}.
\end{align}

We claim that conditioned on the event $\mathcal{A}_n,$ the event $\mc  E^{\ms{uniq}}_n$ occurs with high probability and thus the giant component $\ccn$ is well-defined. To verify  this, we apply  Proposition \ref{cluster in rec}  to the subgraph $\Tilde{B}_n$ of   $\Z^2$ with $V(\Tilde{B}_n ) = V(B(n))$ and $E(\Tilde{B}_n ) = E(  B(n) ) \sm \bar{E}$. To see that
 $\Tilde{B}_n$ satisfies the assumption in Proposition \ref{cluster in rec},
for $i = -k,\cdots,k-1$, define the rectangles  
\begin{align} \label{412}
   R_i&:= [ - n,n]_\Z \times [i \lfloor n/k \rfloor +1,(i+1) \lfloor n/k \rfloor ]_\Z,
\end{align}
where we set $-k\lfloor n/k \rfloor := - n $ and $k\lfloor n/k \rfloor :=  n $, and define
\begin{align*}
  R&:=  [n -\lfloor \tau n \rfloor+1,n]_\Z \times [-n,n]_\Z
\end{align*}
(see Figure \ref{fig} for the illustration).
Then,  $\Tilde{B}_n$  is a union of these rectangles (in both vertex and edge sets sense) which satisfies the overlapping condition in Proposition \ref{cluster in rec}.
Hence, $\P( \mc E_\delta (\Tilde{B}_n)) \ge 0.9$ for   large $n$.

Since $V(\Tilde{B}_n) = V(B(n))$, the giant cluster $\cc$ in Proposition \ref{cluster in rec}  satisfies  $|\cc| \ge \kappa (2n)^2$  for small enough $\de>0$ (recall that $\kappa \in (0,\theta)$ is a constant). In addition, by the property (1) in  Proposition \ref{cluster in rec}, 
 the size of any other cluster in $\Tilde{B}_n$ is at most $(\log n)^5.$   Hence, 
\begin{align} \label{21}
\cA_n \cap \mc E_\delta (\Tilde{B}_n)\Rightarrow \mc E^{\ms{uniq}}_n.
\end{align}
Note that  the events $\cA_n$ and  $\mc E_\delta (\Tilde{B}_n)$ are independent since $\cA_n \in \sigma(\bar{E})$ and $ \mc E_\de(\Tilde{B}_n)  \in \sigma(E(\Tilde{B}_n))$. 

~

Next, we prove  that for large  enough  $k\in \N$ and small  enough $\tau>0$,  the event $\cA_n$ implies that the isoperimetric constant of \emph{any} open circuit is arbitrarily large.
\bepr \label{prop: upper tail lower}
\label{lower prop}
For any constant $t>0$, there exist $k\in \N$ and $\tau>0$
 such that for sufficiently large $n$, under the event $\cA_n \cap \mc E_\delta (\Tilde{B}_n)$,
\begin{align}
 n\wtpn \ge t.
\end{align}
\enpr
\bep
Let $c_0>0$ be the constant from \eqref{372} and define
  $ \zeta:= {c_0^2}/{t^2}.$
Also, let $k\in \N$ and $\tau>0$ be constants such that
  $ 2 c_0\sqrt{k} >t$ and $   16\sqrt{k} \tau < c_0\zeta 
$.
We claim that under the event $\cA_n \cap \mc E (\Tilde{B}_n)$, the isoperimetric constant of any   circuit $\g\in \mathcal C_n$ is large:
\begin{align} \label{421}
	n\frac{|\g|}{|\ccn \cap \vol(\g) |} \ge  n\frac{|\g|}{|\vol(\g) |} \ge t.
\end{align}
By \eqref{333}, recalling    $ \zeta= {c_0^2}/{t^2}$, the above  inequality holds for any {$\g \in \mc C_n$} satisfying $|\vol(\g)| \le \zeta n^2$ .

Hence, it suffices to verify \eqref{421} for any {$\g \in \mc C_n$} such that $|\vol(\g)| \ge \zeta n^2$. 
Set $U_i:= \vol(\g) \cap R_i$ as a subset of $\Z^2$ (see \eqref{412} for the definition of $R_i$). {Since  $ \{U_i\}_{i=-k,\cdots,k-1}$ form a partition of $\vol(\g)$,} 
\begin{align} \label{410}
  |\vol(\g)| = \sum_{i=-k}^{k-1} | U_i|.
\end{align}
Also, setting 
$$L_i:=[n-\lfloor \tau n \rfloor+1,n]_\Z \times \{i \lfloor n/k \rfloor+1, (i+1) \lfloor n/k \rfloor   \},$$ 
{observe that $U_i$ is a vertex-disjoint union of enclosed regions of some circuits $\g_i^{(j)}$s and (possibly) some subsets of  $L_i$ (see the red region in Figure \ref{fig}  for the illustration).} Define $\partial U_i$ to be the union of these circuits $\g_i^{(j)}$s.
Since  $\partial U_i \sm  L_i$s ($i=-k,\cdots,k-1$) are disjoint subsets of $\g$,
\begin{align} \label{411}
  |\g| \ge \sum_{i=-k}^{k-1} |\partial U_i| -{4 k \tau n}.
\end{align}
In addition, since $|L_i| \le 2 \tau n,$
\begin{align*}
 |U_i| \le  \sum_j  |\vol(\g_i^{(j)})| +2\tau n. 
\end{align*}
Furthermore, by the (discrete) isoperimetric inequality \eqref{372}, for each  circuit $\g_i^{(j)}$,
\begin{align*}
     |\g_i^{(j)}| \ge c_0 \sqrt{|\vol(\g_i^{(j)})|} \ge c_0 \frac{\sqrt{k}}{n}   |\vol(\g_i^{(j)})| ,
\end{align*}
where we used $ |\vol(\g_i^{(j)})| \le |R_i| \le n^2/k$ in the last inequality. 
Therefore, 
taking the union over all circuits $\g_i^{(j)}$s in  $\partial U_i$,
\begin{align} \label{34}
 |\partial U_i|  =   \sum_j   |\g_i^{(j)}|  \ge   \sum_j  c_0 \frac{\sqrt{k}}{n}   |\vol(\g_i^{(j)})|  \ge c_0 \frac{\sqrt{k}}{n}   (|U_i|-2\tau n).
\end{align}
Thus, for large enough $n$,
\begin{align*}
   |\vol(\g)| \overset{\eqref{410}}{=}  \sum_{i=-k}^{k-1} | U_i| \overset{\eqref{34}}{\le} \sum_{i=-k}^{k-1} \Big(  \frac{n}{c_0\sqrt{k}}  |\partial U_i| +2\tau n\Big)& \overset{\eqref{411}}{\le}  \frac{n}{c_0\sqrt{k}} ( |\g| + 4k \tau n )+ 4k\tau n \le \frac{n}{c_0\sqrt{k}} |\g| +\frac12|\vol(\g)| ,
\end{align*}
where we used conditions $16\sqrt{k} \tau < c_0\zeta$  and $|\vol(\g)| \ge \zeta n^2$ in the last inequality. Therefore, recalling the condition $2 c_0\sqrt{k} >t$, we obtain \eqref{421}.
\enp

We now conclude the proof of the lower bound for the upper tail.
\bep [Proof of the lower bound in \eqref{theorem4}]
For a given $t>0$, let $k\in \N$ and $\tau>0$ be
 such that Proposition \ref{prop: upper tail lower} holds.
By \eqref{21} and Proposition \ref{prop: upper tail lower},
 \begin{align*}
   \P (n\wtpn \ge t\mid \mc E_n^{\ms{uniq}}) \ge   \P (\{n\wtpn \ge t\}\cap \mc E_n^{\ms{uniq}}) & \ge \P (\cA_n \cap \mc E_\delta (\Tilde{B}_n) ) =\P(\cA_n ) \P (  \mc E_\delta (\Tilde{B}_n)) ,
 \end{align*}
where we used the independence of $\cA_n $ and $ \mc E_\de (\Tilde{B}_n)$  in the last identity.
Since $\P(\cA_n) = e^{-\Theta(n)} $ and $\P ( \mc E_\delta (\Tilde{B}_n)) \ge 0.9$ for large $n$, we conclude the proof.
 \enp

\subsection{Upper bound on the upper tail}
\label{sec:up}
In the present section, we establish the upper bound \eqref{theorem4} in Theorem \ref{thm:main}, i.e.  the upper tail large deviations of the isoperimetric constant.    
From now on, we say that a sequence of the events $\{E_n\}_{n \ge 1}$ occurs \emph{with very high probability} ($\ms{wvhp}$), conditionally on $\mc E^{\ms{uniq}}_n,$ if the probability of the complement   decays exponentially fast,  i.e. $$\limsup_{n \rightarrow \infty }  \frac{1}{n} \log  \P(E_n^c  \mid \mc E^{\ms{uniq}}_n ) < 0.$$

Throughout this section, we condition on the event $\mc E^{\ms{uniq}}_n$ so that the giant cluster $\ccn$ is well-defined.
We prove the large deviation upper bound by verifying that for any $\e>0$,  $\ms{wvhp}$, one can  take $\g_0\in \mathcal C_n$  such that 
$$  n \frac{|\g_0|}{|\ccn \cap \vol(\g_0) |}  \le (1 + \e)  \frac{\xi}{\sqrt{2}\th}.
 $$

For $\e > 0$,  let $\AA_{n, \e}^{\ms W}$ be the event that there exists $\g_0 \in \mathcal C_n$ with $|\vol(\g_0)|\le |B(n)|/2$ such that  
 \begin{align} \label{431}
	 |\ccn \cap \vol(\g_0)|\ge \tfrac12 (1 - 4\e)\th |B(n)|\quad \text{ and }\quad
	 {|\g_0|} \le(1 + 4\e)\th^{-1/2}\xi\sqrt{|\ccn \cap \vol(\g_0)|}.
 \end{align}
 The first condition  is the density property of $\ccn$ in the region $\vol(\g_0)$ and the second condition is about the isoperimetric property of the circuit $\g_0$.
We claim that  $ \AA_{n, \e}^{\ms W}$ occurs $\ms{wvhp}$.

\bepr[Concentration for the Wulff shape]
\label{lem:uu1}
For any   $\e > 0$, the event $ \AA_{n, \e}^{\ms W}$ occurs $\ms{wvhp}$.
\enpr

Assuming Proposition \ref{lem:uu1}, we conclude the proof of the upper bound in \eqref{theorem4}.

\bep[Proof of Theorem \ref{thm:main}, upper bound in \eqref{theorem4}]
For any $\e>0,$
the isoperimetric constant of $\g_0 \in \mathcal C_n$ satisfying the condition \eqref{431} is bounded as
\begin{align*}
	\f{\big|\g_0\big|}{|\ccn \cap \vol(\g_0) |} \le \f1{\sqrt{|\ccn \cap \vol(\g_0)|}} (1 + 4\e) \th^{-1/2}\xi\le \f{1 + 4\e}{\sqrt{1 - 4\e}}\f\xi{\sqrt2 \th n}.
\end{align*}
 Since such $\g_0 \in \mathcal C_n$  exists $\ms{wvhp}$ by Proposition \ref{lem:uu1}, we conclude the proof as $\e > 0$ is arbitrary.
\enp

We aim to establish   Proposition \ref{lem:uu1}. Recall that $\mu$
denotes the time constant (see Proposition \ref{prop: time constant}). For $\e>0$ and  $x,y\in B^{\text{con}}(n)$, let $\G^{(n)}_\e(x,y)$ be the set of \emph{$\e$-optimal path}, i.e., the collection of  open (simple) paths $\g$ in $\ccn$  connecting $[x]$ and $[y]$  (recall that $[x]$ denotes the closest vertex of $x$ in $\ccn$, see Section \ref{sec 3.2})  satisfying
\begin{align*}
    |\g|_1 \le (1 + \e) \mu(y-x).
\end{align*}
A crucial ingredient in  the proof of Proposition \ref{lem:uu1} is  the following proposition, a sharp quantitative concentration result for the closeness of $\e$-optimal paths to the line segment.  
For any  $x,y\in \R^2$,  denote by $[x,y]$ the line segment connecting $x$ and $y$.
Also 
from now on, we simplify the presentation by abusing notation and considering a path $\g $ in $\Z^2$ also as a subset of $\R^2$.

\bepr[Concentration of $\e$-optimal paths]
\label{cor:eop}
Let  $\e,\e' \in (0,1)$ be  constants. Suppose that  $\{x_n\}_{n\ge 1}$ and $\{y_n\}_{n\ge 1}$ are the collection of points such that  $x_n,y_n\in B^\textup{con}(  (1-\e) n)$ and $|y_n-x_n|_1 \ge \e' n.$   Then, $\ms{wvhp}$,
\begin{align*}
     \inf_{\g \in \G^{(n)}_\e( x_n,y_n)}d_{\ms H}(\g, [x_n, y_n]) \le \e |y_n -x_n| _\infty
\end{align*}
(we set the infimum is infinity if $\G^{(n)}_\e(x_n, y_n) = \es$).

\enpr

Proposition \ref{cor:eop} can essentially be considered as a quantitative refinement of \cite[Proposition 3.2 (2)]{blpr} {where the exceptional probability bound $e^{-\Omega( (\log n)^2)}$, much weaker than $e^{-\Omega(n)}$, was obtained for a different notion of distance (called ``right-boundary'' distance) related to the edge boundary  in \eqref{eq:blpr}.}

We postpone the proof of  Proposition \ref{cor:eop} to the end of this section, and elaborate on how to deduce  Proposition \ref{lem:uu1} from this. Let  $W:=W_\mu$  be  the Wulff shape 
w.r.t. the norm $\mu,$ as introduced in \eqref{wulff}, and set $\wh W := |W|^{-1/2}W$ so that 
$|\wh W|=1$. 
 Inspired by \cite[Proposition 4.2]{blpr},  we show that $\ms{wvhp}$, one can choose an open circuit approximating the \emph{Wulff curve}  
 \begin{align} \label{wul}
     \hgp := \pa \wh W
 \end{align}
  (we assume that $\hgp:[0,1]\rightarrow \R^2$ parametrizes $ \pa \wh W$ in a constant speed,  and we abuse the notation  $\hgp$ both for the curve and its trace).
 Note that $ \hgp$ is a Jordan curve since {$\wh W$ is convex}. Also, $\len_\mu (\hgp) = \xi$.

Moreover, general convexity arguments from  \cite[Lemma 5.1(iv)]{blpr} imply that
\begin{align} \label{inside}
    \wh W \su B^\text{con}(1/\sqrt2).
\end{align}
We briefly explain the argument for the sake of completeness. 
Let $r := \max\{r'  > 0: (r', 0 )\in \wh W\}$, which exists since $\wh W$ is closed. Then, by symmetry and  convexity of $\wh W$, we have $B_1(r) \su \wh W$, where  $B_1(r)$ denotes the $\ell^1$-ball of radius $r$ centered at origin. Thus, $2r^2 \le |B_1(r)| \le |\wh W| = 1$, i.e., $r\le 1/\sqrt 2$. The fact $\wh W \su B^\text{con}(r)$ follows by the following reason: If there exists $x:= (x_1, x_2) \in \R^2$ with $x_1 > r$ such that $x\in \wh W$, then $\textbf{n} \cdot x \le \mu(\textbf{n}) / \sqrt{|W|}$ for all $\ell^2$-unit vectors $\textbf{n}$. By symmetry of $\mu,$
$$ \textbf{n} \cdot (2x_1, 0) = \textbf{n} \cdot ((x_1, x_2) + (x_1, -x_2)) \le 2\mu(\textbf{n})/ \sqrt{|W|},$$
 contradicting the choice of $r$.
Hence, we obtain \eqref{inside}.

~

A central conceptual ingredient in the proof of Proposition \ref{lem:uu1} is the notion of polygonal approximation, which has been exploited in the theory of the Wulff shape \cite{acc,dks}.
For $r > 0$, we define the  \emph{$r$-polygonal approximation} $\PP_r(\la)$ of a rectifiable curve $\la:[0,1]\rightarrow  \R^2$, following \cite[Section 4.2]{blpr}. Let $t_0:=0$ and $x_0 := \la(0)$. Then, set $x_k := \la(t_k)$, where we recursively define 
\begin{align*}
    t_k := \inf\big\{t \in (t_{k - 1}, 1]\co |\la(t) - x_{k - 1}|_\infty > r\big\}.
\end{align*}
We stop the construction at $t_J$ when the infimum is taken over an empty set and then define $x_J := \la(1)$. Then, the  $r$-polygonal approximation of $\la$ is defined to be $\PP_r(\la) :=\poly(x_0, \cdots, x_J)$,  a polygonal curve obtained by the concatenation of line segments $[x_{k-1},x_{k}]$ for $k=1,\cdots,J$.

Since  $\len_\rho  ([x,y])  = \rho(y-x)$ for any norm $\rho$, we have 
    $\len_\rho  (\PP_r(\la))  = \sum_{k=1}^J \rho (x_k - x_{k-1})$.
In addition,
\begin{align}\label{ine}
     \len_\rho (\PP_r(\la))  \le \len_\rho (\la) .
\end{align}
Note that  $\PP_r(\la)$ may not be simple even when $\la$ is simple. In order to define the notion of $\intt (\la)$ for a rectifiable closed curve $\la$ which is \emph{not} simple,   let $w_\la(x)$ be the winding
number of $\la$ around $x$, and  then define 
\begin{align} \label{hu}
   \hull(\la) := \la \cup \{x \notin \la: w_\la(x) \text{ is odd}\}.
\end{align}
We present some useful properties of the  polygonal approximation. We define the $\ell^\infty$-Hausdorff metric $\dist$ and  $d_{\ms H}$ as follows:  For compact sets $A,B \subseteq \R^2$ and $x\in \R^2$,
\begin{align*}
    \dist(x,B) := \inf_{y\in B} |x-y|_\infty
\quad \text{ and }\quad
    d_{\ms H}(A,B):= \max \big\{\sup_{x\in A} \dist(x,B)  , \sup_{y\in B}  \dist(y,A) \big\}. 
\end{align*}

\bel [Lemmas 4.3 and  4.4 in \cite{blpr}] \label{poly}
The following properties hold.

1. Let $\lambda$ be a closed rectifiable  curve in $\R^2$. Then, for any $r>0,$
\begin{align} \label{poly1}
    d_{\ms H}(\hull(\la) ,\hull(\PP_r(\la)) )\le r.
\end{align}

2. Let $\la$ be a closed {polygonal} curve in $\R^2$. Then, for any  $\e>0$ and any norm $\rho$ on $\R^2$, there exists a 
\emph{simple} closed polygonal curve $\la'$
such that
\begin{align*}
    |\hull (\la ) \Delta \intt (\la') | < \e \quad  \textup{and}\quad  \len_\rho (\la')\le \len_\rho (\la)+\e,
\end{align*}
where $\Delta$ denotes the symmetric difference.
\enl
From  \eqref{poly1}, we deduce that there exists a constant $C_0>0$ such that for any $r > 0$ and any Jordan curve $\la$ in $\R^2$,
{\begin{align} \label{diff2}
    \big||\intt(\la)| -| \hull(\PP_r(\la))|\big | \le C_0r (\len_{|\cdot|_1}(\la) + r).
\end{align}}

%
%PRF UU1
%

Equipped with this preparation, assuming Proposition \ref{cor:eop}, one can deduce Proposition \ref{lem:uu1}.
\bep[Proof of Proposition \ref{lem:uu1}]
Recall the (normalized) Wulff shape $\hat{W}$ and the Wulff curve $\hat{\g} =  \pa \wh W$ defined in \eqref{wul}.
For $r >0$, denote by $\PP_r(\hgp) = \poly(x_0, \cdots, x_J)$ the $r$-polygonal approximation of  $\hgp$ (we have $x_J=x_0$ since $\hgp$  is closed). Recalling that the  Wulff shape $W$ is the unit ball in the dual norm $\mu^*$ of $\mu$ (see \eqref{dual}), we first claim that  for any $\e>0$, there is $r>0$ small enough such that
\begin{align} \label{55}
   1 - \e^2  \le \sqrt{|W|} \min_{x \in \PP_r(\hgp)}\mu^*(x) \le  1 .
\end{align} 
The second inequality immediately follows from the convexity of $\hat{W}$. To show the first inequality,  if $x \in [x_i, x_{i + 1}]$, then
$$\mu^*(x) \ge \mu^*(x_i) - \mu^*(x_i - x) = 1/ \sqrt{|W|} - |x_i - x|_\ff \mu^*\big((x_i - x)/|x_i - x|_\ff\big).$$
The claim now follows from the observation that $|x_i - x|_\ff \le r$ and $\sup_{u\in \R^2:  |u|_\ff = 1} \mu^*(u) <\infty$.

~

 Denoting $N := (1 -\e)\sqrt2n$,
 we aim to construct  $\g_0 \in \mathcal C_n$ such that 
 \begin{align}
	 \label{eq:wc}
 \max_{x \in \g_0}\dist\big(x,   {N\PP_r(\hat{\g})}\big) \le \e^2 N  \quad \text{ and }\quad |\g_0| \le (1 + \e) N \cdot \len_\mu (\PP_r(\hgp)),
 \end{align}
  where for a curve $\la$ in $\R^2$, $N \la $ denotes a curve obtained by a dilation of $\la$ by $N$.

 We claim that any such $\g_0$ satisfies \eqref{431} $\ms{wvhp}$. To achieve this, we  show that the first property in \eqref{eq:wc} implies (in a deterministic sense)  that for small enough $\e>0$, 
\begin{align}
    \label{eq:oeb}
    (1 - \e)N\hat W \su \vol(\g_0)\su (1 + \e)N\hat W .
\end{align}
Indeed, if the first inclusion does not hold, then recalling that $W$ is the unit ball w.r.t.  $\mu^*$, there would exist $x\in \g_0$ with $\mu^*(x) \le (1 - \e)N / \sqrt{|W|} $. Moreover, by the first property in \eqref{eq:wc}, there exists $x' \in N\PP_r(\hat \g)$ (i.e. $\mu^*(x') \ge   (1 - \e^2) N/ \sqrt{|W|}$, see \eqref{55}) with $|x - x'|_\infty \le \e^2 N$. However, then,
$$(\e - \e^2) N / \sqrt{|W|} \le  \mu^*(x') - \mu^*(x) \le \mu^*(x' - x) \le \e^2 N \sup_{u\in\R^2 : |u|_\infty = 1} \mu^*(u),$$
yielding a contradiction for small $\e>0$.  
Next, to verify the second inclusion in \eqref{eq:oeb}, it suffices to show that  $\max_{x \in \g_0}\mu^*(x) \le (1 + \e)N/ \sqrt{|W|}$. Note that for any $x\in \g_0$, by the first property in \eqref{eq:wc}, there exists $x' \in N\PP_r(\hat \g)$ with $|x - x'|_\infty \le \e^2 N$. Thus, we deduce \eqref{eq:oeb} by noting that
$$\mu^*(x)\le \mu^*(x') + \mu^*(x - x') \overset{\eqref{55}}{\le} N/ \sqrt{|W|} + |x - x'|_\ff \mu^*\big((x - x')/|x - x'|_\ff\big) \le  (1 + \e)N/ \sqrt{|W|}.$$

 Hence, for any $\g_0 \in \mathcal C_n$ satisfying the property \eqref{eq:wc}, by \eqref{lower} in  Proposition \ref{lem:cdens} (conditions $\diam(\g_0) \ge \e^2 n$ and $ n |\g_0|/|\vol(\g_0)| \le  \e^{-2}$ hold by \eqref{eq:wc} and \eqref{eq:oeb}),  for small enough $\e>0$, $\ms{wvhp}$, 
 $$ |\ccn \cap \vol(\g_0)| \ge (1 - \e^2)\th  | \vol(\g_0)|    \overset{\eqref{eq:oeb}}{\ge}  
(1 - \e^2)\theta \cdot  (1-\e)^2N^2 |\hat W| \ge \tfrac12 (1 - 4\e)\th|B(n)|.$$
Also, noting that $ \len_\mu (\PP_r(\hgp))  \le \len_\mu (\hgp)    = \xi$   by  \eqref{ine}, this together with  \eqref{eq:wc} imply that
 $$|\g_0| \le (1 + \e) \xi N \le  (1 + \e) \xi\cdot (1-\e)(1 - 4\e)^{-1/2}\th^{-1/2}\sqrt{|\ccn \cap \vol(\g_0)|} \le (1 + 4\e) \th^{-1/2}\xi\sqrt{|\ccn \cap \vol(\g_0)|},$$
verifying the condition \eqref{431}. Note that
\begin{align*}
    |\vol(\g_0)| \overset{\eqref{eq:oeb}}{\le} (1+\e)^2 N^2 = 2(1+\e)^2 (1-\e)^2 n^2 \le \frac{1}{2}|B(n)|.
\end{align*}

Therefore, it remains to show that $\ms{wvhp}$, there exists $\g_0 \in \mathcal C_n$ satisfying \eqref{eq:wc}. 
 Note that $|Nx_k - Nx_{k-1}|_\infty = \Theta(n)  $ and $Nx_k \in B^\text{con}((1-\e )n)$ (see \eqref{inside} and recall $N = (1-\e) \sqrt{2} n$). Hence, by Proposition \ref{cor:eop} along with a union bound, $\ms{wvhp}$, there exist  open (simple) paths $\g^{(k)}$ connecting $[Nx_{k - 1}]$ and $[Nx_k]$ {in $\ccn$} such that  for all $k =1,\cdots,J$, 
$$d_{\ms H}\big(\g^{(k)}, [Nx_{k - 1}, Nx_k]\big)\le \e^2 N \quad \text{ and }\quad |\g^{(k)}| \le  (1+\e) N\cdot \mu(x_k - x_{k - 1}).$$
 Then, noting the fact {$\len_\mu  ([x,y])  = \mu(y-x)$},  we concatenate the paths $\{\g^{(k)}\}_{ k = 1,\cdots,J}$ as $(\cdots (\g^{(1)} * \g^{(2)})\cdots) * \g^{(J)}$ to yield   $\g_0 \in \mathcal C_n$ 
 satisfying \eqref{eq:wc}. Here, the concatenation  of   simple paths $\g = (u_0, \cdots, u_n)$ and $ \g' = (v_0, \cdots, v_m)$ (assume that $u_n=v_0$)  is defined to be $\g * \g':= (u_0, \cdots, u_k, v_{\ell + 1}, \cdots, v_m)$, where
$k:= \min\{i\co u_i \in \g'\}$ and $\ell$ is a unique index such that $v_\ell = u_k$. Note that the concatenation  $\g * \g'$ is also simple.
\enp

\subsection{Proof of Proposition \ref{cor:eop}}
We conclude this section by establishing  Proposition \ref{cor:eop}. We use the following fact: Since all norms in $\R^2$ are equivalent,  denoting by   $A:= \max_{\textbf{n}\in \R^2, |\textbf{n}|_1 = 1}\mu(\textbf{n})$,
\begin{align} \label{upp}
    |\mu(x) - \mu(y)| \le A |x-y|_1 ,\quad \forall x,y\in \Z^2.
\end{align}
In addition,  for $S_1, S_2\su \R^2$, we write $ S_1+  S_2 :=\{ s_1 +s_2\in \R^2 \co s_1\in  S_1, s_2\in  S_2\}$.
\bep[Proof of Proposition \ref{cor:eop}]
We suppress the notation $n$ in $x_n$ and $y_n$ for the sake of readability. 
We construct a specific $\e$-optimal path $\g$ in $\ccn$ from  $[x]$ to $[y]$ which is close to the line segment $[x, y]$.
To achieve this goal, following the idea in \cite[Proposition 3.2]{blpr}, for $M\in \N$, define  $u_k := (1 - kM^{-1})x + kM^{-1}y \in \R^2$ ($k=0,\cdots,M$) to be  the  points equi-partitioning the segment $[x,y]$. We show that if $M$ is sufficiently large, then  the geodesic $\g^{(k)}$ from $[u_{k - 1}]$ to $[u_k]$ ($k=1,\cdots,M$), which is a priori \emph{not} necessarily in $\ccn$ (i.e., $\g^{(k)}$ can leave  $B(n)$), is close to the line segment $[x,y]$ and thus contained in $\ccn$. In addition, we prove that 
the concatenation
$\g := (\cdots(\g^{(1)} * \g^{(2)})* \cdots) *\g^{(M)}$  
satisfies the desired properties. Note that $\g$ is simple, since each geodesic  $\g^{(k)}$ is simple.
\medskip

\textbf{Step 1. $\g$ is contained in $\ccn$ and $d_{\ms H}(\g, [x, y]) \le \e |y-x|_\infty$.}   Recall that $x,y\in B^\textup{con}(  (1-\e) n)$ and $|y-x|_1 \ge \e' n.$  This ensures that for every $k=1,\cdots, M,$ $u_{k-1},u_k \in B^\textup{con}(  (1-\e) n) $ and $|u_k - u_{k-1}|_1 \ge \e' n/M.$ Hence,
 by \eqref{311} and \eqref{310},  $\ms{wvhp}$,
\begin{align} \label{bulk}
 [u_{k-1}],[u_{k}] \in  [x, y] + B^{\text{con}}(\e \e' n/4)
\end{align}
and  
\begin{align*}
    |[u_k] -[u_{k-1}]|_1  = \Theta(n).
\end{align*}
The chemical distance  $D([u_{k-1}],[u_k])$ can be estimated using \eqref{31200}: $\ms{wvhp}$, 
\begin{align} \label{312}
| D([u_{k-1}],[u_k]) - \mu([u_{k}] - [u_{k-1}]) | \le  \e  \mu([u_{k}] - [u_{k-1}])/3.
\end{align}
Also,
by \eqref{upp}, 
$$\big|\mu([u_{k}] - [u_{k-1}]) - \mu(u_{k} - u_{k-1})\big| \le A \big|([u_{k}] - [u_{k-1}]) - (u_{k} - u_{k-1})\big|_1 .$$
The inequality \eqref{310} along with the fact $ \mu(z) \ge |z|_1 $ imply that $\ms{wvhp}$, the above quantity is bounded by  $\e \mu(u_{k} - u_{k-1}) /3.$
Hence, this together with 
\eqref{312} give that $\ms{wvhp}$, 
\begin{align} \label{contain}
D([u_{k-1}],[u_{k}]) \le  (1+\e) \mu(u_{k} - u_{k-1}) .
\end{align}
We claim that  for large enough $M$, under the events \eqref{bulk} and \eqref{contain}, 
\begin{align} \label{contain1}
    \g^{(k)} \text{ is contained in } [x,y]+B^{\text{con}}(\e |y-x|_\infty).
\end{align} 
Suppose that there exists  a lattice point $w$ in  $ \g^{(k)}$ such that $w \notin [x,y]+B^{\text{con}}(\e |y-x|_\infty)$. Since $|y-x|_\infty \ge |y-x|_1/2 \ge \e' n/2$,
% we have \red{$d_{\ms H}(w, [x, y]) >  \e |y-x|_\infty/2 \ge \e \e' n/4 $.} Hence, 
this together with \eqref{bulk}  would imply that
\begin{align} \label{notcontain}
    D([u_{k-1}],[u_{k}])  =|\g^{(k)}|  \ge |[u_{k-1}]-w|_1 + |w-[u_{k}]|_1  >  \e  \e' n/2.
\end{align}
Note that if  $M > 16A/(\e \e')$, then
 \begin{align*}
     \mu(u_{k} - u_{k-1}) \le A|u_{k} - u_{k-1}|_1  = A|x-y|_1 / M \le 4An/M < \e \e' n/4,
 \end{align*}
which together with \eqref{contain} and \eqref{notcontain} yield a contradiction for small  $\e>0$. Hence, we deduce  \eqref{contain1}.

Therefore, by  \eqref{contain1} 
applied to all $k=1,\cdots,M$, $\ms{wvhp}$,
\begin{align} \label{113}
    \g\text{ is contained in }[x,y]+B^{\text{con}}(\e |y-x|_\infty).
\end{align}
Noting that {$[x,y]+B^{\text{con}}(\e |y-x|_\infty) \subseteq B(n)$} (recall that $x,y\in B^\textup{con}((1-\e)n)$), $\g$ is entirely contained in $\ccn.$  In addition,
since $|x-[x]|, |y-[y]| \le \e|y-x|_\infty/4$ $\ms{wvhp}$ (by \eqref{311}) and $ \g$ is a path connecting $[x]$ and $[y]$, from \eqref{113}, we obtain
\begin{align*}
    [x,y] \text{ is contained in } \g + B^{\text{con}}(\e |y-x|_\infty). 
\end{align*}
This together with \eqref{113} imply that $\ms{wvhp}$, $d_{\ms H}(\g, [x, y]) \le \e |y-x|_\infty$.  

~

\textbf{Step 2. $\bs {\g \in \G^{(n)}_\e(x,y)}$.} 
By \eqref{contain}, $\ms{wvhp}$,
\begin{align*}
|\g|_1  \le  \sum_{k =1}^M |\g^{(k)}|_1 = \sum_{k =1}^M D([u_{k-1}],[u_{k}]) \le (1+ \e )\sum_{k =1}^M \mu(u_k - u_{k - 1}) = (1+  \e)\mu(y - x) ,
\end{align*}
where in the last identity, we used the fact that $u_1,\cdots,u_M$ lie on the segment $[x,y].$  
\enp

 \section{Lower tail large deviations}
 \label{sec:low}

 \subsection{Lower bound on the lower tail}
In this section, we  lower bound 
 the lower tail probability of the isoperimetric constant, by establishing the lower bounds in \eqref{theorem2} and \eqref{theorem3}.
First, the lower bound in \eqref{theorem2} immediately follows from the observation that if all edges in $B(n)$ are open (which happens with probability $e^{-\Theta( n^2)}$), then the event $\mc E^{\ms{uniq}}_n$ holds and $\ccn=B(n)$. Thus, by taking  $\g\in \mathcal{C}_n$ parametrizing $ \partial B( \lfloor n/\sqrt{2}\rfloor )$ in the variational problem \eqref{eq:isop},
\begin{align*}
    \wtpn \le
 \f{4\sqrt{2} n}{|B(\lfloor n/\sqrt{2}\rfloor )|} 
 \le \frac{2\sqrt{2}}{n} .
\end{align*}

Next, we prove  the lower bound in \eqref{theorem3}, by showing that  for any $\e>0,$
\begin{align} \label{501}
 \P\Big(n\wtpn \le \frac{2\sqrt{2}}{\theta}+\e \mid \mc E^{\ms{uniq}}_n \Big)  \ge e^{-\Theta(n)}.
\end{align}
Setting  $m:=\lfloor n/\sqrt{2}\rfloor$,  let us denote by $F$ the collection of edges in $B(n)$ having at least one endpoint in $\partial B(m) .$ 
Then,  define the event
\begin{align*}
   \mc{H}:= \{\text{All edges in} \
   F\ \text{are open}\}.
\end{align*}
Then, $\P(\mc{H}) = e^{-\Theta(n)}.$

Let $\delta>0$ be a small constant which will be chosen later.
Let $\cE_1$  be the event that  there exists a crossing cluster $\cc_1\su B(m-1)$ with $|\cc_1| \ge (\theta-\de) |B(m-1)|$ and the diameter of any other cluster  in $B(m-1)$ is at most $(\log n)^2.$
Then  by Theorem \ref{cluster in square} applied to $B(m-1)$, we have $\P(\cE_1 ) \ge 0.9$.  
In addition, by considering the annulus  $B(n) \sm B(m) $  as a union of four ``long'' rectangles (see Figure \ref{fig2}), let $\cE_2$ be the event that there exists a strongly crossing cluster $\cc_2\su B(n) \sm B(m) $ with $|\cc_2| \ge (\theta-\de)|B(n) \sm B(m)|$ and the diameter of any other cluster in 
 $B(n) \sm B(m)$ is at most $(\log n)^2.$ Then,
by Proposition \ref{cluster in rec} applied to the annulus  $B(n) \sm B(m) $, we have $\P(\cE_2 ) \ge 0.9$.

\usetikzlibrary{patterns}
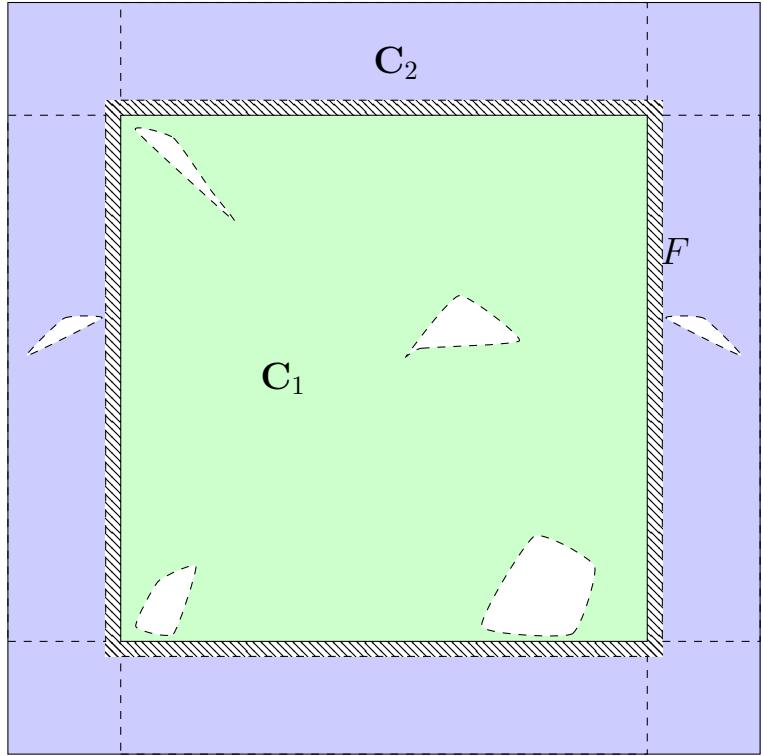
\begin{figure}
	\centering
 \begin{tikzpicture}

	 \fill[blue!20] (-5,-5) rectangle (5,5);
	 \fill[white] (-3.5,-3.5) rectangle (3.5,3.5);
	 \fill[green!20] (-3.5,-3.5) rectangle (3.5,3.5);

	   \fill[white] plot[smooth cycle, tension = 0.3] coordinates{(-3.3, -3.3) (-2.8, -3.4) (-2.5, -2.5)  (-3.0, -2.7)};
	   \draw[dashed] plot[smooth cycle, tension = 0.3] coordinates{(-3.3, -3.3) (-2.8, -3.4) (-2.5, -2.5)  (-3.0, -2.7)};

\fill[white] plot[smooth cycle, tension = 0.3] coordinates{(1.3, -3.3) (2.5, -3.4) (2.8, -2.5)  (2.0, -2.1)};
\draw[dashed] plot[smooth cycle, tension = 0.3] coordinates{(1.3, -3.3) (2.5, -3.4) (2.8, -2.5)  (2.0, -2.1)};

	 \fill[white] plot[smooth cycle, tension = 0.3] coordinates{(0.3, 0.3) (0.5, 0.4) (1.8, 0.5)  (1.0, 1.1)};
	 \draw[dashed] plot[smooth cycle, tension = 0.3] coordinates{(0.3, 0.3) (0.5, 0.4) (1.8, 0.5)  (1.0, 1.1)};

\fill[white] plot[smooth cycle, tension = 0.3] coordinates{(-3.3, 3.3) (-2.8, 3.2) (-2.3, 2.5)  (-2.0, 2.1)};
\draw[dashed] plot[smooth cycle, tension = 0.3] coordinates{(-3.3, 3.3) (-2.8, 3.2) (-2.3, 2.5)  (-2.0, 2.1)};

	 \fill[white] plot[smooth cycle, tension = 0.3] coordinates{(-4.75, 0.3) (-4.25, 0.8) (-3.75, 0.8) };
	 \draw[dashed] plot[smooth cycle, tension = 0.3] coordinates{(-4.75, 0.3) (-4.25, 0.8) (-3.75, 0.8) };

	 \fill[white] plot[smooth cycle, tension = 0.3] coordinates{(4.75, 0.3) (4.25, 0.8) (3.75, 0.8) };
	 \draw[dashed] plot[smooth cycle, tension = 0.3] coordinates{(4.75, 0.3) (4.25, 0.8) (3.75, 0.8) };

\fill[white] (-3.7,-3.7) rectangle (3.7,-3.5);
\fill[white] (-3.7,3.5) rectangle (3.7,3.7);
\fill[white] (-3.7,-3.7) rectangle (-3.5,3.7);
\fill[white] (3.5,-3.7) rectangle (3.7,3.7);

\draw[dashed] (-3.5,3.7) rectangle (3.5,5);
\draw[dashed] (-3.5,-3.7) rectangle (3.5,-5);
\draw[dashed] (3.7,-3.5) rectangle (5,3.5);
\draw[dashed] (-3.7,-3.5) rectangle (-5,3.5);

\fill[pattern=north west lines, pattern color=black] (-3.7,-3.7) rectangle (3.7,-3.5);
\fill[pattern=north west lines, pattern color=black] (-3.7,3.5) rectangle (3.7,3.7);
\fill[pattern=north west lines, pattern color=black] (-3.7,-3.7) rectangle (-3.5,3.7);
\fill[pattern=north west lines, pattern color=black] (3.5,-3.7) rectangle (3.7,3.7);

	   \coordinate[label=-180:{\Large{$\cc_1$}}] (A) at (-0.9, 0);
	   \coordinate[label=-180:{\Large{$F$}}] (A) at (4.2, 1.7);
	   \coordinate[label=-180:{\Large{$\cc_2$}}] (A) at (0.6, 4.2);
	 \draw (-5,-5) rectangle (5,5);
	 \draw (-3.5,-3.5) rectangle (3.5,3.5);

 \end{tikzpicture}
	\caption{Giant clusters $\cc_1$ (green), $\cc_2$ (blue) and region $F$ (hatched). Smaller clusters are enclosed by dashed curves.}
 \label{fig2}
\end{figure}

According to the crossing property of $\cc_1$ and $\cc_2,$ under the event $\mc{H}\cap \cE_1 \cap \cE_2$, $F\cup \cc_1 \cup \cc_2$ (i.e., the graph consisting of edges in $F$, $\cc_1$, and $\cc_2$) becomes a cluster in $B(n)$ which satisfies
\begin{align} \label{515}
  |F\cup \cc_1 \cup \cc_2|  \ge (\theta-2\de) (2n)^2 > \kappa (2n)^2
\end{align}
(recall that $\kappa \in (0,\theta)$ is a constant).
Also, by the  uniqueness of the giant component condition imposed on $\cE_1$ and $\cE_2$, the event $\mc E^{\ms{uniq}}_n $ holds, and thus we have $\ccn = F \cup  \cc_1 \cup \cc_2.$

Hence, under the event $\mc{H}\cap \cE_1 \cap \cE_2$,  by taking $\g \in \mathcal C_n$ parametrizing $\partial B(m) $ in \eqref{eq:isop} (which satisfies $ |\vol(\g)|\le \frac12 |B(n)|$), for small enough $\de>0,$
\begin{align*}
    \wtpn \le  \f{8m}{| \ccn \cap B(m)|} \le 
 \f{8m}{|\cc_1|} \le \frac2{\theta -\de} \frac1{m} \le \frac{1}{n} \Big(\frac{2\sqrt{2} }{\theta}+\e\Big).
\end{align*}
Note that by the   independence of $\cE_1, \cE_2$ and $ \mathcal{H}$,  
\begin{align}
  \P(\mc{H}\cap  \cE_1\cap \cE_2 ) = \P(\mc{H} ) \P( \cE_1) \P( \cE_2 ) \ge  (0.9)^2 e^{-\Theta(n)}.
\end{align}
Therefore,
we obtain \eqref{501}.

%
%SEC LUP
%
\subsection{Upper bound on the lower tail}
\label{sec:lup}
In this section, we establish \eqref{theorem1} and  the upper bounds in \eqref{theorem2} and \eqref{theorem3}, i.e., the upper bound on the lower tail large deviations of the isoperimetric profile.

\subsubsection{Proof of \eqref{theorem1}}
 The statement \eqref{theorem1} is a consequence of the deterministic fact that for any $\de>0$, for sufficiently large $n$ and any circuit $\g$ in $\Z^2$ such that $|\vol(\g)| \le |B(n)|/2$,
\begin{align} \label{522}
  \f{n|\g|}{|\vol(\g)|} \ge 2\sqrt{2} - \de.
\end{align}
To see this, by Lemma \ref{discrete iso}, there exists $R>1$ such that for any circuit $\g$ in $\Z^2$ such that $|\vol(\g)| > R$ and $|\g| < |\vol(\g)|^{2/3}$ satisfies $|\g| \ge (4-\de) \sqrt{|\vol(\g)|}$. This implies
\begin{align} \label{525}
  \f{n|\g|}{|\vol(\g)|} = \frac{n}{\sqrt{|\vol(\g)|}} \cdot \f{|\g|}{ \sqrt{|\vol(\g)|}} \ge \frac{n}{\sqrt{|B(n)|/2}}  \cdot (4-\de)  \ge  2\sqrt{2} -\de.
\end{align}
In addition, if $\g$ satisfies either $|\vol(\g)|\le R$ or $|\g|  \ge  |\vol(\g)|^{2/3}$, then for sufficiently large $n$,
\begin{align*}
 \f{n|\g|}{|\vol(\g)|} \ge \frac{n}{R} \ge 2\sqrt{2}- \de
\end{align*}
(in the  former case $|\vol(\g)|\le R$) or
\begin{align*}
\f{n|\g|}{|\vol(\g)|} \ge \frac{n}{|\vol(\g)|^{1/3}} \ge 2\sqrt{2}-\de
\end{align*}
(in the latter case $|\g| \ge |\vol(\g)|^{2/3}$). Therefore, we conclude the proof of \eqref{522}.

\subsubsection{Proof of the upper bound in \eqref{theorem2}}
We show that for any $\e>0$,
\begin{align} \label{5210}
   \P\Big(\inf_{ \substack{\g \in \mathcal C_n\\ |\vol(\g)| \le |B(n)|/2 } } \f{n|\g|}{|\ccn \cap \vol(\g) |} \ge \frac{2\sqrt{2}}{\theta}-\e\mid \mc E^{\ms{uniq}}_n \Big) \ge 1-e^{-\Omega(n^2)}.
\end{align} Throughout the proof, we assume the occurrence of the event $\mc E^{\ms{uniq}}_n$ so that the cluster $\ccn$ is well-defined.  Since we aim to lower bound the quantity $n|\g|/|\ccn \cap \vol(\g)|$  uniformly in   $\g \in \mathcal C_n$, it suffices only to consider circuits $\g \in \mathcal C_n$ such that  for some small   constant $\de>0$, 
\begin{align} \label{assume1}
    n|\g|/|\vol(\g)|  \le \de^{-1}.
\end{align}
 If $\diam(\g) \le \de n$, which implies $ |\vol (\g)| \le  2(\diam(\g))^2 \le 2\de^2 n^2$,  then by \eqref{333}, $n|\g|/|\ccn \cap \vol(\g)| \ge n|\g|/|\vol(\g)| \ge c_0/  
(\sqrt{2} \de)$. Thus,  we henceforth also assume that  for some small constant $\de>0$, 
 \begin{align} \label{assume2}
     \diam(\g) \ge \de n.
 \end{align} 
  
Then, by Proposition \ref{lem:cdens}, with probability at least $1-e^{-\Theta(n^2)}$, for any  $\g \in \mathcal C_n$ with $|\vol(\g)| \le |B(n)|/2$ satisfying \eqref{assume1} and \eqref{assume2},
\begin{align} \label{523}
  \f{n|\g|}{|\ccn \cap \vol(\g) |} &=\f{n|\g|}{|\vol(\g)|}\cdot   \frac{|\vol(\g)|}{|\ccn \cap \vol(\g) |} \overset{\eqref{522}}{\ge} (2\sqrt{2}-\delta) \cdot \frac{1}{\theta+\de} .
\end{align}
Since $\delta>0$ is arbitrary, we conclude the proof of \eqref{5210}.

\subsubsection{Proof of the upper bound in \eqref{theorem3}}

The key strategy 
is to approximate the open circuit $\g \in \mathcal C_n$ by a suitable Jordan curve $\la$ in $\R^2$ such that $\ms{wvhp}$, $\len_\mu(\la)$ is essentially at most $|\g|$. Once this is achieved, the  isoperimetric inequality applied to the Jordan curve $\la$ concludes the proof.  Throughout this section, we condition on the event $\mc E^{\ms{uniq}}_n$ so that $\ccn$ is well-defined.   

\bel 
\label{lem:pr41}
 For any  small enough $\de>0$,  the following event, which we call  $\aac$, occurs $\ms{wvhp}$. For any $\g\in \mathcal C_n$ satisfying $ \diam(\g) \ge \de n$ and  $ n|\g|/|\vol(\g)|  \le \de^{-1}$, there exists a Jordan curve $\la$ in $\R^2$ such that
\begin{align} \label{511}
    \big||\vol(\g)| - | \intt(\la)|\big| \le \de^{}{|\vol(\g)|}  \ \textup{ and } \  |\g| \ge (1 - \de^2)\len_\mu(\la).
\end{align}
\enl
Assuming Lemma \ref{lem:pr41}, we verify the upper bound in \eqref{theorem3}.

\bep[Proof of the upper bound in \eqref{theorem3}]
We show that for any $\delta>0$,
\begin{align} \label{521}
   \P\Big(\inf_{\substack{\g \in \mathcal C_n\\ |\vol(\g)| \le |B(n)|/2}} \f{n|\g|}{|\ccn \cap \vol(\g) |} \ge \frac{\xi}{\sqrt{2}\theta}-\delta   \ \mid \mc E^{\ms{uniq}}_n \Big) \ge 1-e^{-\Omega(n)}.
\end{align}
Similarly as in the proof of \eqref{5210}, we only consider circuits $\g\in \mathcal C_n$ with $|\vol(\g)| \le |B(n)|/2$  satisfying \eqref{assume1} and \eqref{assume2}.
Under the event $\{s^+_{n, \de}\le  \de \} \cap \aac ,$ for any such $\g\in \mathcal C_n$, there exists  a Jordan curve $\la$  in $\R^2$ satisfying \eqref{511} (see  Lemma \ref{lem:pr41}), and thus
\begin{align*}
    \f{n|\g|}{ |\ccn \cap \vol(\g)  |} & =   \frac{n}{\sqrt{|\vol(\g)|}} \cdot \f{|\g|}{\sqrt{|\vol(\g)|}}\cdot \frac{|\vol(\g)|}{|\ccn \cap \vol(\g)  |} \\
    &\ge  \frac{n}{\sqrt{|B(n)|/2}} \cdot   \f{(1 - \de^2)\len_\mu(\la)}{\sqrt{(1 - \de)^{-1}|\intt(\la)|}}  \cdot  \frac{1}{\theta + \delta} \ge  \frac{\xi}{\sqrt{2}\theta}-c(\delta)
\end{align*}
for some $c(\delta)>0$ with $\lim_{\delta \downarrow 0} c(\de)=0$,
where the last inequality follows from the isoperimetric inequality $\len_\mu(\la) \ge  \xi \sqrt{|\intt(\la)|}$ (see \eqref{iso percolation}). 

 Since $\{s^+_{n, \de}\le  \de \} \cap \aac $ occurs  $\ms{wvhp}$ by Proposition \ref{lem:cdens} and Lemma \ref{lem:pr41}, we conclude the proof.
\enp

It remains to prove Lemma \ref{lem:pr41}.

\bep[Proof of Lemma \ref{lem:pr41}]
Setting $r :=  \lfloor \de^3 n \rfloor,$ define the event  
\begin{align*}
    \mathcal{G}_{n, \de}:= \{D(x,y) \ge (1 - \de^2/2)\mu(y - x),\quad \forall x, y\in {\ccn}  \text{ with } |y - x|_1\ge r \}.
\end{align*}
{By  Proposition \ref{ld distance} together with   translation invariance and a union bound,} 
the event $\mathcal{G}_{n, \de}$  occurs $\ms{wvhp}$.

We show that under the event $\mathcal{G}_{n, \de}$, for  any $\g\in \mathcal C_n$ satisfying $ \diam(\g) \ge \de n$  and $ n|\g|/|\vol(\g)|  \le \de^{-1}$,  there exists a  Jordan curve  $\la$  in $\R^2$ with the asserted properties \eqref{511}. Regarding $\g$ as a Jordan curve in $\R^2$ (parametrized by a constant speed), let $\la = \poly(z_0, \cdots, z_J)$ be an $r$-polygonal approximation of $\g$. As $r \in \N$ and  $\g\in \mathcal C_n$, we have $z_k \in \ccn$ for  $k=0,\cdots,J$. 

Since $\la$ may not be simple, we  first show that $\la$ satisfies
\begin{align} \label{999}
    \big||\vol(\g)| - | \hull(\la)|\big| \le \de {|\vol(\g)|} /2 \ \textup{ and } \  |\g| \ge (1 - \de^2)\len_\mu(\la)
\end{align}
(see  \eqref{hu} for the definition of $\hull(\la)$). We conclude the proof once this is achieved,  as any closed  polygonal curve can be approximated by a closed ``simple''    polygonal curve (see the second statement in Lemma \ref{poly}). 
 
First,   we have  
$$\big||\vol(\g)| - |\hull(\la)|\big| \overset{\eqref{diff}}{\le}  \big||\intt(\g)| - |\hull(\la)|\big| + |\g|  \overset{\eqref{diff2}}{\le}{C_0 r ( |\g| + r )  + |\g| }.$$
Since $ n|\g|/|\vol(\g)|  \le \de^{-1}$, we have
 $r|\g| \le  \de^{3}n|\g|    \le \de^2|\vol(\g)|.$
 Also, using the fact $|\g| \ge  \diam(\g) \ge \de n,$ we have  $r^2 \le \de^6 n^2 \le \de^5  n |\g| \le \de^4  |\vol(\g)|.$
 Hence, for small $\de>0$, we  verify the first condition in \eqref{999}.

Next, since $z_k\in \ccn$ and  
$|z_k - z_{k - 1} |_1 \ge |z_k - z_{k - 1} |_\infty = r$ for $k=1,\cdots,J$, under the event $\mathcal{G}_{n, \de}$, 
$$|\g| \ge \sum_{k=1}^{J} D(z_{k - 1}, z_k) \ge (1 - \de^2/2) \sum_{k=1}^{J} \mu(z_k - z_{k - 1})= (1 - \de^2/2)(\len_\mu(\la) - \mu(z_J - z_{J - 1})).$$
We lower bound the above quantity. Note that using the fact $|z_k - z_{k - 1} |_1 \le 2|z_k - z_{k - 1} |_\infty \le 2r$ for all $k=1,\cdots,J$ (indeed, the last inequality becomes equality for $k=1,\cdots,J-1$) along with a triangle inequality, we deduce  $\len_{|\cdot|_1}(\la) \ge \diam(\g)  - 4r$. Hence, for some constant $A'> 0$, 
\begin{align*}
    \len_\mu(\la) \ge A' \len_{|\cdot|_1}(\la) \ge A' (\diam(\g)  - 4r) \ge A' \de n/2,
\end{align*} 
where the last inequality follows from the condition $\diam(\g) \ge \de n$. Thus,
    $$\mu(z_J - z_{J - 1}) \overset{\eqref{upp}}{\le} A  |z_J - z_{J - 1}|_1 \le 2A  r \le 2 A \de^3n \le C \de^2 \len_\mu(\la).$$
Therefore, for small $\de>0$, we deduce
the second condition in \eqref{999}.
\enp

%
%PPENDIX
%
\section{Appendix}
\label{sec:App}
In the Appendix, we prove Proposition \ref{cluster in rec},  Lemma \ref{lemma 3.4} and  Proposition \ref{prop: time constant dominance}. 

\subsection{Proof of  Proposition \ref{cluster in rec}} 
In this section, we establish Proposition \ref{cluster in rec}, using the idea in \cite[Proposition 2]{pp}. Indeed,
we take a suitable collection of ``dominos'' (of size $(\log n)^2 \times 2(\log n)^2$  and  $2(\log n)^2 \times (\log n)^2$) covering $R_n$. 
Such suitable choice of dominos  allows the crossing clusters in each domino to be connected together, and controls the diameter of any other cluster  in $R_n$.
As Proposition \ref{cluster in rec} aims to construct a cluster which further satisfies a density condition (2), we need an additional ingredient, namely Lemma \ref{domino}  below, which ensures the existence of the crossing cluster in the rectangle (of side-length ratio $1:2$ and $2:1$) whose density is essentially at least $\theta$.

\bel \label{domino}
Let $\de>0$ be a constant and $D_n$ be a rectangle of size $2n\times n$ (or $n \times 2n$). Define $\mc E_\de (D_n)$ to be the event such that
 there exists a crossing cluster $\cc$ in $   D_n$ satisfying
 $ |\cc| / 2n^2 \ge \theta - \de .$
Then, there exists a constant $c>0$ such that for sufficiently large $n$,
\begin{align*}
  \P (\mc E_\de (D_n))  \ge 1-e^{-cn}.
\end{align*}

\enl

Assuming Lemma  \ref{domino},  we  establish Proposition \ref{cluster in rec}.

\bep [Proof of Proposition \ref{cluster in rec}]
Define $m:= \lfloor (\log n)^2 \rfloor.$
A rectangle $D$ is called a \emph{domino} if it is of size either $m\times 2m$  or $2m \times m$. Let us take dominos of the  following particular form  contained in $R_n$:
\begin{align*}
  [im, (i+2)m] _\Z\times [jm, (j+1)m]_\Z,\quad [im, (i+1)m] _\Z\times [jm, (j+2)m]_\Z.
\end{align*}
Since these dominos may not cover the region close to the boundary of $R_n$, we  take additional dominos (with the smallest  cardinality) which cover the whole region $R_n$.
Let $\{D_\ell\}_{\ell \in L}$ be the collection of such dominos. Then, for some constant $C>0,$
\begin{align} \label{bound}
    |L| \le C\frac{n^2}{(\log n)^4} .
\end{align}
By Lemma \ref{domino}, for each $\ell\in L$, there exists an event $\mc E_{\de/2}(D_\ell)$ satisfying
\begin{align} \label{321}
  \P(\mc E_{\de/2}(D_\ell)) \ge 1-{e^{-c (\log n)^2}}
\end{align}
under which there is a crossing cluster $\cc(D_\ell)\su D_\ell$ such that
\begin{align} \label{322}
  \frac{|\cc(D_\ell)|}{2m^2} \ge \theta - \frac\de 2.
\end{align}

By the second condition imposed on the structure of $R_n,$ under the intersection of events $\mc E_{\de/2}(D_\ell)$ for all $\ell\in L$, which we call $\mc{F}$, clusters $\cc(D_\ell)$s are connected together, yielding a   cluster $\cc\su   R_n$.
By \eqref{bound} and  \eqref{321}, for sufficiently  large $n$,
{\begin{align}
   \P(\mc{F}) \ge 1- C\frac{n^2}{(\log n)^4} \cdot  e^{-c(\log n)^2} \ge 1- e^{-c(\log n)^2/2} .
\end{align}}

We verify that the cluster $\cc$ satisfies the property (2). First, the strongly crossing property of $\cc$ is a direct consequence of the crossing property of $\cc(D_\ell)$s. Next, to deduce the bound \eqref{12}, take a maximal collection of ``disjoint'' dominos $D_\ell$s  in $R_n$, whose union is denoted by $R_n'$. Then, by \eqref{322},   for large $n,$
\begin{align*}
  |\cc| \ge \big(\theta - \tfrac\de 2\big) |R_n'| \ge (\theta- \de) |R_n|,
\end{align*}
where the last inequality follows from the condition that the size of each rectangle $R_n^\ell$ is $\Theta(n) \times \Theta(n)$ and the size of domino is poly-log in $n$.

Finally, the property (1) is a consequence of the crossing property of each cluster $\cc(D_\ell) \su D_\ell$.

\enp

We now provide the proof of Lemma \ref{domino}. 
\begin{proof}[Proof of Lemma \ref{domino}]
Assume that $D_n = [0,2n]_\Z \times [0,n]_\Z,$ and define the  squares
\begin{align*}
  B_1:=[0,n]_\Z \times [0,n]_\Z,\quad B_2:=[\lfloor n/2 \rfloor,\lfloor 3n/2 \rfloor]_\Z \times [0,n]_\Z,\quad B_3:=[n,2n]_\Z \times [0,n]_\Z.
\end{align*}
   By Theorems 4 and 5 in \cite{pp} with $\phi_n := \de n$, with probability at least $1-e^{-cn}$,   there exists a crossing cluster $\cc_i$ in $   B_i$ ($i=1,2,3$) such that 
   \begin{align} \label{331}
 \frac{|\cc_i|}{n^2} \ge\theta-\de.
   \end{align}
   These clusters $\cc_i$s ($i=1,2,3$) can be connected together to form a crossing cluster $\cc$ in $  D_n$, provided that there is a vertical crossing in the rectangles $[\lfloor n/2 \rfloor,n]_\Z \times [0,n]_\Z$ and $[n,\lfloor 3n/2 \rfloor]_\Z \times [0,n]_\Z$ respectively. Such vertical crossings exist with probability at least $1-e^{-c'n}$ (by a duality and the exponential tail decay of the radius of an open cluster in the subcritical regime \cite[Theorem 5.4]{grim}). Under such an event, by \eqref{331}, the cluster $\cc$ satisfies \begin{align*}
     |\cc| \ge |\cc_1|+|\cc_3| \ge (\theta-\de) \cdot 2n^2,
   \end{align*}
   where the first inequality follows from the disjointness of $\cc_1$ and $\cc_3$.
\end{proof}

\subsection{Proof of Lemma \ref{lemma 3.4}}

In this section, we prove Lemma \ref{lemma 3.4} by following the argument in \cite[Theorem 2]{ds} where the density estimate was obtained for the infinite cluster $\ccf$ (see  \eqref{ld}).
 
\begin{proof}[Proof of Lemma \ref{lemma 3.4}]

   Let  $K>0$ be a large constant  which will be chosen later, and suppose that the side lengths of $R$ are at least $M:= 2K/\sqrt{\theta}$. Then, for any cluster $\cc\su   \Z^2$ such that $|\cc \cap R| \ge (\theta+\de) |R|$,
   \begin{align} \label{343}
     |\cc| \ge(\theta+\de) |R| \ge (\theta+\de) M^2 \ge 4K^2.
   \end{align}
Assume that each side length of $R$ is multiple of $K$, and let $D_1,\cdots,D_L$ be $(K \times K)$-squares defining a partition of $R$.
  For $x\in R,$ let $\ell_x$ be the unique $\ell \in \{1,2,\cdots,L\}$ such that $x\in D_\ell$, and define
  \begin{align*}
    \zeta(x) := \begin{cases}
1 \quad \text{$x$ is connected to $\partial D_{\ell_x}$ through open edges in $D_{\ell_x}$} ,\\
0 \quad \text{otherwise}.
    \end{cases}
  \end{align*}
Then, for any cluster $\cc$, $\zeta(x)=1$ for any $x\in \cc \cap R$. This is because  $|\cc| \ge 4K^2 >|D_{\ell_x}| $ (see \eqref{343}), which implies that $\cc$ contains a point outside  $D_{\ell_x}$. Hence,
\begin{align} \label{341}
  \sum_{x\in R} \zeta(x)\ge |\cc \cap R| .
\end{align}
Motivated by the fact that
the random fields $(\zeta|_{D_\ell})_{\ell=1,\cdots,L}$  are i.i.d., we write
\begin{align} \label{342}
  \frac1{|R| } \sum_{x\in R} \zeta(x) = \frac1 L \sum_{\ell =1}^L \Big(\frac1{|D_\ell|} \sum_{x\in D_\ell} \zeta(x)\Big).
\end{align}
We upper bound the one-point expectation of $\zeta$. By  translation invariance, for any $\ell=1,\cdots,L$ and $x\in D_\ell $ such that $\dist(x,\partial D_\ell) \ge \sqrt K,$ 
\begin{align*}
	\P(\zeta(x)=1) \le \P(|\cc(o)|\ge \sqrt K ) =: q_K,
\end{align*}
	where $\cc(o)$ denotes the cluster containing the origin $o$.
Thus,
\begin{align*}
  \E \Big[\frac1{|D_\ell|} \sum_{x\in D_\ell} \zeta(x)\Big]\le q_K + \frac{4 K \sqrt{K}}{K^2} = q_K + \frac4{\sqrt K }. 
\end{align*}
Since $\lim_{K \rightarrow  \infty} q_K = \theta,$  for any  $\de>0,$ the above quantity is less than $\theta+\de$ for a large enough constant $K$.
Hence, for such  $K$, by the large deviation estimate for the sum of i.i.d. (bounded) random variables, there exists  a constant {$c = c(\de)>0$} such that
\begin{align*}
  \P\Big( \frac1 L \sum_{\ell =1}^L \Big(\frac1{|D_\ell|} \sum_{x\in D_\ell} \zeta(x)\Big) \ge \theta+\de \Big) \le e^{-cL} .
\end{align*}
Therefore, recalling the bound \eqref{341} and the identity \eqref{342}, we deduce \eqref{340}.

 In the case when some side length of $R$ is not a multiple of $K,$ take the maximum size of rectangles $R'\su R$ whose side lengths are multiple of $K$.  Then, there exists $M = M(K,\de)>0$ such that if the side lengths of $R$ are at least $M$, the event $|\cc \cap R| \ge(\theta+\de)|R|$ implies
\begin{align*}
  |\cc\cap R'| \ge |\cc \cap R| - |R \sm R'| \ge (\theta+\de)|R| - \frac\de 2 |R|  \ge \Big (\theta+ \frac\de 2 \Big) |R|.
\end{align*}
Here, the second inequality follows from the maximality of $R'$ and the fact that side-lengths of $R$ are sufficiently large.
By the above argument applied to $R'$ and noting that $|R| \le 2|R'|$ (for sufficiently large $M$), we conclude the proof.
\end{proof}

%
%SS TIMDOM
%
\subsection{Proof of Proposition \ref{prop: time constant dominance}} \label{sec 7.3}
We conclude the Appendix by establishing Proposition \ref{prop: time constant dominance}.  A crucial ingredient is the following proposition which  provides the exponential  probability bound on the event that the chemical distance, in the ``almost'' horizontal and vertical directions, is comparable to the $\ell^1$-distance.
\bepr \label{prop: time constant inequality}
Consider supercritical percolation in $\Z^2.$
Then, there exist constants $C,c, \tau,\zeta>0 $ such that for any $x = (x_1,x_2)\in \Z^2$ such that $|x_2| \le \tau |x_1| $ or $ |x_2| \ge \frac1{\tau}|x_1| $,
\begin{align*}
  \mathbb{P} (D(o,x) \le (1+\zeta)|x|_1) \le Ce^{-c |x|_1}.
\end{align*}

\enpr

Given Proposition \ref{prop: time constant inequality}, one can conclude the proof of Proposition \ref{prop: time constant dominance}.

\bep[Proof of Proposition \ref{prop: time constant dominance}]
Let $C,c, \tau,\zeta>0 $ be constants from Proposition \ref{prop: time constant inequality}. Let $x\in \Z^2$ be such that $|x_2| \le \tau |x_1| $ or $ |x_2| \ge \frac1{\tau}|x_1| $, and let $(T_{n,x})_{n\ge 1}$ be the increasing sequence of positive integers such that $\{T_{n,x} x \longleftrightarrow \infty\}$.
 Then, by  Proposition \ref{prop: time constant inequality},
\begin{align*}
  \mathbb{P} (D(o,T_{n,x} x) \le (1+\zeta)|T_{n,x}x|_1) \le Ce^{-c |T_{n,x}x|_1}.
\end{align*}
By the Borel-Cantelli lemma, almost surely, for sufficiently large $n$,
\begin{align*}  
  D(o,T_{n,x} x) > (1+\zeta)|T_{n,x}x|_1.
\end{align*}
This together with Proposition  \ref{prop: time constant} imply that  $\mu(x) \ge (1+\zeta) |x|_1.$
\enp

Finally, we prove Proposition \ref{prop: time constant inequality}. A similar version of Proposition \ref{prop: time constant inequality} was previously obtained in the context of FPP as follows.

\bepr[Proposition 3.1 in \cite{mar}] \label{prop: fpp}
Consider the FPP in $\Z^2$ with i.i.d.~edge distributions having a common law $\nu$ supported on $[0,\infty)$. Assume that $\inf \textup{supp} (\nu)= 1$ and $1/2 \le \nu(\{1\}) <1$.  For $x,y\in \Z^2,$ let $D^{\textup{FPP},\nu}(x,y)$ be the infimum of all passage times from $x$ to $y$.
Then, there exist constants $C,c,\tau,\zeta>0 $ such that for any $x= (x_1,x_2)\in \Z^2$ such that $|x_2| \le \tau |x_1| $ or $ |x_2| \ge \frac1{\tau}|x_1| $,
\begin{align} \label{fpp estimate}
  \mathbb{P} (D^{\textup{FPP},\nu}(o,x) \le (1+\zeta)|x|_1) \le Ce^{-c|x|_1}.
\end{align}
\enpr

The chemical distance in supercritical percolation can  be   regarded as FPP with i.i.d.~edge distributions $\nu$ on $[0,\infty) \cup \{\infty\}$ given by $\nu(\{1\}) = p$ and $\nu ( \{ \infty\}) = 1-p$ (i.e. the open and closed edge correspond to passage times 1 and $\infty$, respectively).
Although Proposition \ref{prop: fpp} is not directly applicable since  $\nu$ is supported on $\{1,\infty\},$ one can detour this problem by a simple coupling argument. In fact, we consider the FPP with i.i.d.~edge distributions $\nu'$ given by $\nu'(\{1\}) = p$ and  $\nu'(\{2\}) = 1-p$. Since $1/2\le p<1$, by Proposition \ref{prop: fpp}, we have \eqref{fpp estimate} for the first passage time $D^{\textup{FPP},\nu'}$. Since $\nu$ stochastically dominates $\nu'$, there is a coupling between passage times such that $D^{\textup{FPP},\nu} \ge D^{\textup{FPP},\nu'}.$ This concludes the proof of  Proposition \ref{prop: time constant inequality}.

\bibliographystyle{plain}
\bibliography{iso}

\end{document}